\documentclass{article}

\usepackage[english]{babel}

\usepackage[letterpaper,top=2cm,bottom=2cm,left=3cm,right=3cm,marginparwidth=1.75cm]{geometry}

\usepackage{amsmath,empheq}
\usepackage{amsfonts}
\usepackage[dvipsnames]{xcolor}
\usepackage{graphicx}
\usepackage{bm}
\usepackage{amsmath}
\usepackage{subfig}
\usepackage[percent]{overpic}
\usepackage{pict2e}
\usepackage{xspace}
\usepackage{colortbl}
\usepackage[colorlinks=true, allcolors=blue]{hyperref}
\usepackage[siunitx]{circuitikz}
\usepackage[maxbibnames=99]{biblatex}

\usepackage{listofitems}

\usepackage{tikz}
\usetikzlibrary{calc,fit,backgrounds}
\usetikzlibrary{shapes.geometric, arrows}
\usetikzlibrary{chains,
                fit,
                positioning,
                shapes}

\usepackage{color}

\tikzstyle{startstop} = [rectangle, rounded corners, minimum width=3cm, minimum height=1cm,text centered, draw=black, fill=red!30]
\tikzstyle{io} = [trapezium, trapezium left angle=70, trapezium right angle=110, minimum width=2cm, minimum height=1cm, text centered, text width=2cm, draw=black, fill=cyan!20!white]
\tikzstyle{io1} = [trapezium, trapezium left angle=70, trapezium right angle=110, minimum width=3cm, minimum height=1cm, text centered, text width=3cm, draw=black, fill=yellow!40!white]
\tikzstyle{process} = [rectangle, minimum width=1cm, minimum height=1cm, text centered, text width=2.5cm, draw=black, fill=orange!20!white]
\tikzstyle{decision} = [diamond, minimum width=3cm, minimum height=1cm, text centered,text width=1.5cm, draw=black, fill=green!60!white]
\tikzstyle{arrow} = [thick,->,>=stealth]
\tikzstyle{matheq} = [node distance=8.75cm, text width=21em, minimum width=1cm, minimum height=2em, text centered]

\tikzset{>=latex} 
\colorlet{myred}{brown!80!black}
\colorlet{myblue}{teal!80!black}
\colorlet{mygreen}{violet!60!black}
\colorlet{mydarkred}{myred!40!black}
\colorlet{mydarkblue}{myblue!40!black}
\colorlet{mydarkgreen}{mygreen!40!black}

\tikzstyle{node}=[very thick,circle,draw=myblue,minimum size=22,inner sep=0.5,outer sep=0.6]
\tikzstyle{connect}=[->,thick,mydarkblue,shorten >=1]
\tikzset{ 
  node 1/.style={node,mydarkgreen,draw=mygreen,fill=mygreen!25},
  node 2/.style={node,mydarkblue,draw=myblue,fill=myblue!20},
  node 3/.style={node,mydarkred,draw=myred,fill=myred!20},
}
\def\nstyle{int(\lay<\Nnodlen?min(2,\lay):3)}

\newcommand{\openfoam}{Open\nolinebreak\hspace{-.2em}{\color{blue}\Large$\nabla$}\nolinebreak\hspace{-.2em}FOAM\textsuperscript{\textregistered}\xspace}

\title{A hybrid Reduced Order Model to enforce outflow pressure boundary conditions in computational haemodynamics}
\author{Pierfrancesco Siena$^{*}$, Pasquale Claudio Africa$^{*}$, Michele Girfoglio$^{\dag}$, Gianluigi Rozza$^{*}$}

\date{}

\begin{document}
\maketitle

\begin{center}
\small $^{*}$SISSA, International School for Advanced Studies, Mathematics Area, mathLab, via Bonomea 265, I-34136 Trieste, Italy \\
\small $^{\dag}$University of Palermo, Department of Engineering,  Viale delle Scienze, Ed. 8, 90128 Palermo, Italy 
\end{center}

\begin{abstract}
This paper deals with the development of a Reduced-Order Model (ROM) to investigate hemodynamics in cardiovascular applications. It employs the use of Proper Orthogonal Decomposition (POD) for the computation of the basis functions and the Galerkin projection for the computation of the reduced coefficients. The main novelty of this work lies in the extension of the lifting function method, which typically is adopted for treating non-homogeneous inlet velocity boundary conditions, to the handling of non-homogeneous outlet boundary conditions for the pressure, representing a very delicate point in the numerical simulations of the cardiovascular system. Moreover, 
we incorporate a properly trained neural network in the ROM framework to approximate the mapping from the time parameter to the outflow pressure, which in the most general case is not available in closed form. We define our approach as ``hybrid", because it merges physics-based elements with data-driven ones.  
At full order level, a Finite Volume method is employed for the discretization of the unsteady Navier-Stokes equations while a two-element Windkessel model is adopted to enforce a valuable estimation of the outflow pressure. 
Numerical results, firstly related to a 2D idealized blood vessel and then to a 3D patient-specific aortic arch, demonstrate that our ROM is able to accurately approximate the FOM with a significant reduction in the computational cost. 
\end{abstract}

\section{Introduction}

The computational investigation of haemodynamics in the cardiovascular system has been the subject of extensive research. Many studies employ high-fidelity Full-Order Models (FOMs) based on finite element or finite volume  approximation of the governing Partial Differential Equations (PDEs) to capture the complex blood flow dynamics. Although these models provide accurate and detailed insight, they are computationally expensive and time-consuming, especially in parametric settings where multiple simulations are required for varying physical and/or geometrical parameters \cite{africa2024lifex}.
On the other hand, Reduced Order Models (ROMs) have emerged as an efficient alternative to FOMs, offering a significant reduction in the computational cost while maintaining a reasonable level of accuracy \cite{gerner2012certified,patera2007reduced,huynh2012certified,benner2020model,noor1981recent,quarteroni2007numerical,manzoni2012computational,deparis2009reduced}. The crucial element of ROMs is that many problems have an intrinsic dimension significantly lower than the number of degrees of freedom in the discretized system. This reduction in dimensionality is achieved through a two-phase process, namely an \emph{offline phase} and an \emph{online phase}. In the \emph{offline phase}, a database of solutions is computed by solving the original FOM for several parameter values. These solutions are then combined and compressed to construct a reduced space. Then the \emph{online phase} adopts the compressed information from the offline phase to predict reduced solutions for the values of the new parameters, significantly reducing the cost of the computation. 

ROMs include two main categories: data-driven approaches (see, e.g., \cite{wang2019non,hesthaven2018non,guo2019data,xiao2015non,audouze2013nonintrusive}) and equations-based methods (see, e.g., \cite{quarteroni2011certified,quarteroni2015reduced,hesthaven2016certified}). In the former family, the surrogate approximations are derived exclusively from the available data, leveraging interpolation or regression techniques. These methods suffer from a complete lack of error estimation theory and the accuracy is significantly influenced by the quantity and quality of the input data. On the contrary, a strong theoretical background is available for the latter family \cite{quarteroni2011certified}, for which the surrogate approximation is achieved by solving a reduced-dimensional version of the original system of PDEs. A well-known subset of equations-based methods is represented by the Galerkin ROM \cite{rowley2004model,iollo2000stability,rozza2007stability}, where the projection of the PDEs system on the reduced space is performed, providing a differential-algebraic system of equations for the reduced coefficients that can be solved using standard iterative methods. 

In this paper, we explore the application of ROMs to the investigation of haemodynamics in cardiovascular applications. We employ the unsteady Navier-Stokes equations coupled with a three-element Windkessel model to enforce realistic outflow pressure conditions \cite{westerhof2009arterial,girfoglio2021non,Girfoglio2020}. A Finite Volume (FV) method is adopted for the space discretization due to its widespread adoption in commercial and open-source codes, as it is robust and exhibits local conservation properties.

We consider time as the only parameter. Parametric studies will be
addressed in future work. 
Many works have already analyzed the characteristics of data-driven ROM techniques to study the blood dynamics in the coronary arteries, aortic arch, and cardiac chambers (see, i.e., \cite{siena2023data,siena2023fast,sienanew2024,balzotti2023reduced,balzottidata2022,Girfoglio2020,buoso2019reduced}), reporting a significant computational speedup with respect to the reference FOM. 
However, in the context of data-driven methods, the main challenge is represented by the large amount of data required for an accurate reconstruction. 
For this reason, the ROM adopted in this work is based on the Galerkin projection of the FOM onto the reduced space, preserving the physics of the problem. 
The reduced space is computed using the Proper Orthogonal Decomposition (POD) method \cite{chatterjee2000introduction}. 
 
Within this framework, our focus is on the control of the outflow boundary condition for the pressure. This is because in cardiovascular applications, the enforcement of this condition
is essential to achieve realistic numerical simulations \cite{fevola2021optimal,vignon2010outflow}. 
When employing a projection-based ROM (as in our case), nonhomogeneous boundary conditions are typically not preserved at the reduced-order level \cite{saddam2017pod}. 
Moreover, they are not explicitly included in the ROM and, therefore, cannot be directly controlled. In order to address this issue, two main approaches can be found in the literature: the penalty method \cite{lorenzi2016pod,kalashnikova2010stability,sirisup2005stability} and the lifting function or control function method \cite{saddam2017pod,graham1999optimal,gunzburger2007reduced}. In order to avoid the sensitivity analysis needed to choose the penalty parameter in the former technique, the lifting function approach is adopted here. 
The lifting function method aims to homogenize the boundary conditions of the basis functions in the reduced space. In particular, the idea is to find a homogeneous reduced basis space so that at the reduced level any boundary condition can be employed. That goal could be achieved by subtracting the nonhomogeneous boundary conditions of the problem from the FOM snapshots collected during the offline phase; then, the boundary value is added again in the reconstruction of the reduced solution. 
Although this method is widely used and applied for the nonhomogeneous inflow boundary conditions for the velocity field \cite{saddam2017pod}
, to the best of our knowledge, its application to the pressure field is explored here for the first time. 

A further complication faced in this work is as follows. 
Unlike the inlet velocity where the boundary value (or the functional form in the case of nonuniform boundary condition) is given, 
for the outflow pressure, we have only discrete values obtained from the FOM snapshots which are collected at a certain frequency. 
So, we have introduced a properly trained feedforward neural network to map from time instances to the output pressure values when the ROM is run for time step values different from the FOM one. 

The rest of the paper is structured as follows. In Section \ref{sec:fom}  the FOM is introduced, reporting all the details regarding the FV discretization. Section \ref{sec:rom} presents our ROM approach. The POD-Galerkin method combined with the lifting function is introduced, followed by a brief explanation of feedforward neural networks. The numerical results are reported in Section \ref{sec:results} both for an idealized 2D blood vessel and for a patient-specific 3D aortic arch. Finally, in Section \ref{sec:conclusion}, we draw conclusions and discuss future perspectives of the current study.
\section{The full order model}
\label{sec:fom}

\subsection{Problem formulation}
\label{sec:navier-stokes}

Consider the motion of an incompressible, viscous, and Newtonian fluid governed by the unsteady Navier-Stokes equations within the spatial domain $\Omega \in \mathbb{R} ^3$ over the time interval $(t_0, T]$: 

   \begin{empheq}[left=\empheqlbrace]{alignat=4}
       \frac{\partial \bm{u}(\bm{x}, t)}{\partial t} + \nabla \cdot \left(\bm{u}(\bm{x}, t) \otimes \bm{u}(\bm{x}, t)\right) - \nabla \cdot (\nu \nabla \bm{u}(\bm{x}, t)) +\nabla p(\bm{x}, t) & =  0 & \quad \mbox{ on } \Omega \times (t_0, T], \label{eq:NS1} \\
      \nabla \cdot \bm{u}(\bm{x}, t) & =  0 & \quad \mbox{ on } \Omega \times (t_0, T], \label{eq:NS2}
   \end{empheq}
with $\bm{u}=\bm{u}(\bm{x}, t)$ the velocity vector, $p=p(\bm{x}, t)$ the kinematic pressure, and $\nu$ the kinematic viscosity. 
In the remainder of this paper, to simplify the notation, we will omit the spatial-temporal dependence of the variables. 

Let us introduce the Reynolds number $Re$  defined as:
\begin{equation}\label{eq:Re}
    Re = \frac{U L}{\nu},
\end{equation}
where $U$ and $L$ are characteristic macroscopic velocity and length for the flow field at hand. In the applications discussed in this paper, the Reynolds number $Re$ is small enough so that the flow can be considered laminar; therefore, no additional modeling is required. 



Let $\Gamma_{i}$, $\Gamma_{o}$, and $\Gamma_{w}$ be the inlet, the outlet, and the wall of the domain $\Omega$, respectively, such that $\Gamma_i \cup \Gamma_w \cup \Gamma_o = \partial\Omega$ and $\Gamma_i \cap \Gamma_w = \Gamma_i \cap \Gamma_o = \Gamma_w \cap \Gamma_o = \emptyset$. 
For the velocity field a nonhomogeneous Dirichlet boundary condition is enforced on the inlet section:
\begin{equation}
    \label{eq:non-hom_dirich}
    \bm{u} = \bm{u}_D \quad \mbox{ on } \Gamma_{i} \times (t_0, T],
\end{equation}
a homogeneous Neumann boundary condition is applied on the outlet section:
\begin{equation}
    \nabla \bm{u} \cdot \bm{n} = 0 \quad \mbox{ on } \Gamma_{o} \times (t_0, T],
\end{equation}
and finally a no-slip boundary condition is employed on the wall
\begin{equation}
    \bm{u} = 0 \quad \mbox{ on } \Gamma_{w} \times (t_0, T].
\end{equation}

\begin{figure}
\begin{center}
\begin{circuitikz}[american]
\draw[] (-0.5,0.5) -- (-0.50,-0.5);
\draw [-{Stealth[scale=1.5]}] (-0.5,0) -- (0.3,0);
\draw[color=purple] (0,0) to[R, l^=$R_{p}$] (3,0);
\draw[color=yellow!70!black] (3,0) to[C, l=$C$,-] (3,-2.5);
\draw[color=green!50!black] (3,0) to[R, l^=$R_{d}$] (6,0);
\draw[] (0,0.3) node[] {$Q$};
\draw (-0.8,0) node[] {$p$};
\draw (3,0.2) node[] {$p_{p}$};
\draw (3,-2.7) node[] {$p_{d}$};
\draw (6.3,0) node[] {$p_{d}$};
\draw[] (-0.45,-0.8) node[] {$\Gamma_o$};
\end{circuitikz}
\end{center}
\caption{Sketch of the three-element Windkessel model.}
\label{fig:windkessel}
\end{figure}
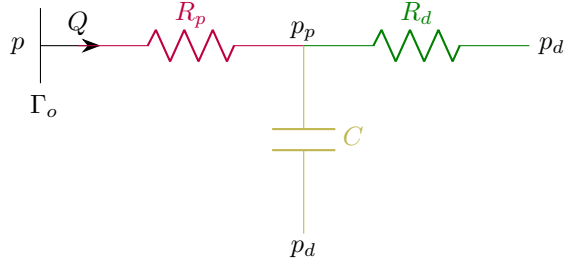

On the other hand, for the pressure field, a homogeneous Neumann boundary condition is enforced on the inlet and on the wall sections:
\begin{equation}
    \nabla p \cdot \bm{n} = 0 \quad \mbox{ on } (\Gamma_{w} \cup \Gamma_{i}) \times (t_0, T].
\end{equation}
In order to enforce realistic outflow boundary conditions, at each outlet of the domain, a three-element Windkessel model is employed
\cite{westerhof2009arterial}. This model is composed of a compliance $C$ and two resistances, the proximal resistance $R_p$ and the distal resistance $R_d$. On a specific outlet section, the downstream pressure $p$ is given by the following differential-algebraic equations system: 
\begin{equation}
   \left \{
   \begin{alignedat}{3}
    C \frac{dp_p}{dt} + \frac{p_p - p_d}{R_d}  & = Q & \quad & \mbox{ on } \Gamma_o,\\
    p - p_p & = R_p Q & \quad & \mbox{ on } \Gamma_o,
   \end{alignedat}
   \right .
   \label{eq:Windkessel-FOM}
\end{equation}
where $p_p$ is the proximal pressure, $p_d$ is the distal pressure (assumed to be, as often in literature \cite{nichols2022mcdonald}, $p_d=0$, as it serves as a reference value) and $Q$ is the flow rate through the outlet section. A scheme is shown in Figure \ref{fig:windkessel}.

\subsection{The discrete formulation}\label{sec:discretization}
For space discretization, we adopt the FV method. The domain $\Omega$ is divided into $N_h$ non-overlapping subdomains $\Omega_i$ with $i\in 1, \dots, N_h$. 

By integrating Equation~\eqref{eq:NS1} over each control volume $\Omega_i$ we get:
\begin{equation}
    \int_{\Omega_i} \frac{\partial \bm{u}}{\partial t} d\Omega +  \int_{\Omega_i} \nabla \cdot (\bm{u} \otimes \bm{u}) d\Omega - \nu \int_{\Omega_i} \Delta \bm{u} d\Omega + \int_{\Omega_i} \nabla p d\Omega = 0.
\end{equation}
Then the application of the Gauss-divergence theorem yields:
\begin{equation}
    \int_{\Omega_i} \frac{\partial \bm{u}}{\partial t} d\Omega +  \int_{\partial\Omega_i} (\bm{u} \otimes \bm{u}) \cdot d\bm{A} - \nu \int_{\partial\Omega_i} \nabla \bm{u} \cdot d\bm{A} + \int_{\partial\Omega_i} p d\bm{A} = 0.
\end{equation}
The boundary terms are approximated by means of the following sums over each face $j$ of the cell $\Omega_i$: 

\begin{align}
    & \int_{\partial\Omega_i} p d\bm{A} \approx \sum_{j} p_{i,j} \bm{A}_{ij},  \\
    & \int_{\partial\Omega_i} (\bm{u} \otimes \bm{u}) \cdot d\bm{A} 
    \approx \sum_{j} (\bm{u}_{i,j} \cdot \bm{A}_{ij}) \bm{u}_{i,j} = \sum_{j} \phi_{i,j} \bm{u}_{i,j},  \\
    & \int_{\partial\Omega_i} \bm{\nabla u} \cdot d\bm{A} \approx \sum_{j} (\bm{\nabla u}_i)_j \cdot \bm{A}_{ij}, 
\end{align}
where $p_{i,j}$, $\bm{u}_{i,j}$, and $(\bm{\nabla u}_i)_j$ respectively denote the values of the pressure, the velocity, and the gradient of the velocity at the centroid of the face $j$, $\bm{A}_{ij} = A_{ij} \bm{n}_{ij}$ where $A_{ij}$ is the area of the face $j$ and $\bm{n}_{ij}$ is the normal unit vector to the face $j$, and $\phi_{i,j} = \bm{u}_{i,j} \cdot \bm{A}_{ij}$ is the convective flux associated with the velocity 
through face $j$ of the control volume $\Omega_i$. All these quantities are obtained by means of a second-order accurate scheme consisting of a linear interpolation of the values from cell centers to face centers. 

Now we are going to comment on the time discretization.
Consider an equispaced grid with time step $\Delta t = {(T-t_0)}/{N_T}$, where $N_T$ is the number of time subintervals. At each time $t_n=t_0+n\Delta t$ with $n=1,\dots,N_T$, the time derivative of $\bm{u}$ is approximated with a first-order Euler scheme: 
\begin{equation}
    \frac{\partial \bm{u}}{\partial t} \approx \frac{\bm{u}^{n+1}-\bm{u}^n}{\Delta t},
\end{equation}
where $\bm{u}^n$ is the approximation of the velocity at the time $t_n$. If $\bm{u}_i^{n+1}$ and $\bm{b}_i^{n+1}$ indicate the average velocity and the average source term in the control volume $\Omega_i$,  $\bm{u}_{i,j}^{n+1}$ and $p_{i,j}^{n+1}$ are the velocity and pressure associated with the centroid of face $j$ normalized by the volume of $\Omega_i$, the discretized form of problem  \eqref{eq:NS1} is given by:
\begin{equation}
\frac{1}{\Delta t} \bm{u}_i^{n+1}+\sum_j \phi_{i,j}^n \bm{u}_{i,j}^{n+1}-\nu\sum_j(\nabla \bm{u}_i^{n+1})_j\cdot \bm{A}_{ij}+\sum_j p_{i,j}^{n+1}\bm{A}_{ij}=\bm{b}_i^{n+1}.
\label{FV_form}
\end{equation}
Now we consider  
the semi-discretized form of \eqref{FV_form}, i.e. with the pressure term in continuous form while all the other terms are in discrete form. By using the divergence-free constraint \eqref{eq:NS2} and the Gauss-divergence theorem we get: 
\begin{equation}
\sum_k (\nabla p^{n+1}_i)_k\cdot \bm{A}_{ik} = \sum_k  \Big( -\sum_j \phi_{i,j} \bm{u}_{i,j}^{n+1} + \nu \sum_j(\nabla \bm{u}_i^{n+1})_j\cdot \bm{A}_{ij} + \bm{b}_i^{n+1} \Big)_k \cdot \bm{A}_{ik}.
\label{FV_form1}
\end{equation}
For further details the reader is referred to \cite{girfoglio2019finite,siena2023fast}. 
In order to treat the velocity-pressure coupling in \eqref{FV_form}-\eqref{FV_form1}, the Pressure Implicit with Splitting of Operators (PISO) algorithm \cite{issa1986solution} is adopted. 

The Euler scheme is also used to discretize the Windkessel system \cite{Girfoglio2020}: 
\begin{equation}
\label{wind_discr}
    \left \{
   \begin{alignedat}{3}
    C \frac{p_p^{n+1} - p_p^n}{\Delta t} + \frac{p_p^{n+1} }{R_d} & = Q^n &  \quad & \mbox{ on } \Gamma_o,\\
    p^{n+1} - p_p^{n+1} & = R_p Q^n & \quad & \mbox{ on } \Gamma_o.
   \end{alignedat}
   \right .
\end{equation}




\section{The reduced order model}\label{sec:rom}
In this work we employ a POD-Galerkin ROM: see, e.g., \cite{sirovich1987turbulence,lorenzi2016pod}. Nonhomogeneous boundary conditions are integrated in the reduction step through the lifting function or control function method \cite{fick2018stabilized,graham1999optimal}. This approach is widely used for the velocity field \cite{saddam2017pod}; however, it remains unexplored for the pressure. This is because usually, in classic computational fluid dynamics, homogeneous Dirichlet or Neumann boundary conditions are enforced for the pressure field. On the other hand, in the context of cardiovascular flows, one needs to integrate nonhomogeneous outflow boundary  conditions for the pressure. 
One substantial difference between the employment of the lifting function for the velocity and for the pressure is that, while the velocity is imposed directly (see equation \eqref{eq:non-hom_dirich}), the outflow pressure is the solution of the Windkessel model (see equation \eqref{eq:Windkessel-FOM}); therefore, in general, we lack a closed-form expression for that. This is a problem if we are going to reconstruct the solution for time instances different from the training ones.
In order to overcome this issue, once the FOM has been solved for certain time instances (and the corresponding discrete set of outflow pressure values is obtained by the Windkessel model \eqref{wind_discr}), a properly trained feedforward neural network is introduced in order to get a map from the parameter space to the outflow pressure.  

The most significant assumption in the ROM context is to approximate the solution as a linear combination of basis, $\bm\phi(\bm x)$ and $\psi(\bm x)$, 
depending only on space, and coefficients $a_i(t)$ and $b_i(t)$, 
related only with time: 
\begin{equation}
    \bm u \approx \bm u_{\text{rb}} = \sum_{i=1}^{N_{\bm u}} a_i(t) \bm\phi_i(\bm x),
    \label{u_rb}
\end{equation}
\begin{equation}
    p \approx p_{\text{rb}} = \sum_{i=1}^{N_p} b_i(t) \psi_i(\bm x),
    \label{p_rb}
\end{equation}
where $N_p$ and $N_u$ are the dimensions of the reduced basis space for pressure and velocity, respectively. 
This assumption allows us to decouple the computations into an expensive phase (offline) to be performed only once and a cheap phase (online) to be performed for every new evaluation. We note that in this work we limit ourselves to building a ROM for the time reconstruction of the system, i.e. the time is the only parameter; however, this procedure can be, in principle, extended to any other parametric study of physical and/or geometrical types. 

\begin{figure}
\begin{center}
       \begin{tikzpicture}[node distance=2cm]
        \node (input) [io] {Create the FOM database};
        \node (pod) [process, below of=input] {Compute lifting functions};
        \node (input01) [io, below of=pod] {Homogeneize FOM snapshots};
        \node (input1) [process, below of=input01] {Perform POD/Galerkin projection};
        \node (int) [io, below of= input1] {Interpolate with NN the outflow pressure};
        
        
        \node (online1) [io1, right of=input1, xshift=4cm,yshift=0.8cm] {Compute the ROM coefficients};
        \node (online2) [io1, right of=pod, xshift=4cm,yshift=-0.5cm] {Construct the ROM solution};
        
         \begin{pgfonlayer}{background}
             \node (fit1)[draw, dashed, rounded corners, fill=none, fit=(input.top right corner)(pod)(input01.bottom left corner)(int), inner sep=10] {}; 
             \node (fit1Label)[above=0cm of fit1] {\textbf{Offline / Training phase}};
             
             \node (fit2)[draw, dashed, rounded corners, fill=none, fit=(online2.top right corner)(online1.bottom left corner), inner sep=10] {}; 
             \node (fit1Label)[above=0cm of fit2] {\textbf{Online / Evaluation phase}};
         \end{pgfonlayer}
         
        \draw [arrow] (input) -- (pod);
        \draw [arrow] (pod) -- (input01);
        \draw [arrow] (input01) -- (input1);
        \draw [arrow] (input1) -- (int);
        \draw [arrow] (online1) -- (online2);
        
        \node (p0) at (0.0,-10.0) {}; 
        \node (p1) at (0.0,-10.0) {}; 
        \node (p2) at (6.05,-10.0) {}; 
        \draw [->] (fit1.south) -- (p0.north) -- (p1.north) -- (p2.north) -- (fit2.south);
        
        \node (p3) at (-3.5,2.3) {}; 
        \node (p4) at (0,2.3) {}; 
        \node (p5) at (0,1.5) {}; 

 
    \end{tikzpicture}
\end{center}
    \caption{Schematic flowchart of the ROM framework. 
    }
    \label{fig:flowchart}
\end{figure}
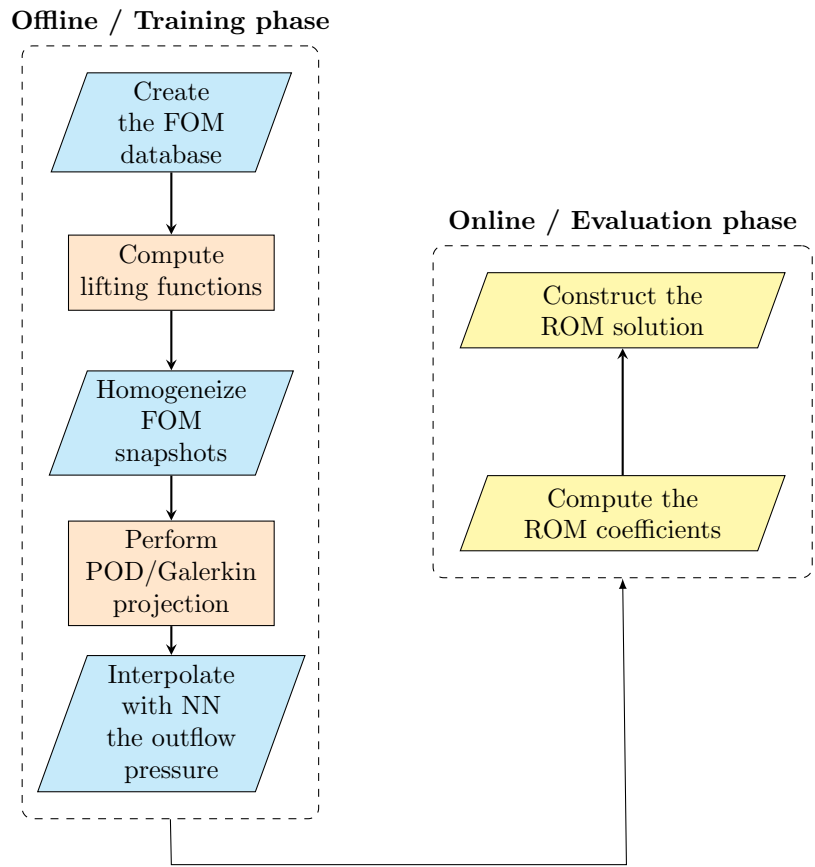

The essential ingredients of the offline and online phases are reported in the following strategy:
\begin{itemize}
    \item {\bf Offline}: Given a discrete set of time instances, the high fidelity solutions, i.e. the FOM solutions, are computed. Firstly, the lifting functions for velocity and pressure are computed and the snapshots are accordingly homogenized (see Sec.~\ref{sec:lift}). The reduced basis space is constructed by the POD algorithm applied to the homogenized solutions (see Sec.~\ref{sec:pod}). 
    Then a dynamical system for the reduced coefficients is obtained through a Galerkin projection of the governing equations \eqref{eq:NS1} and \eqref{eq:NS2} onto the POD modes (see Sec.~\ref{sec:galerkin}). Finally, a neural network is properly trained for mapping the outflow pressure (see Sec. \ref{sec:nn}). 
    \item {\bf Online}: Given a set of new time instances, i.e. different from the training values, the corresponding modal coefficients are computed by solving the dynamical system obtained in the offline phase. Then the ROM solution is obtained 
    by equations \eqref{u_rb} and \eqref{p_rb}. 
\end{itemize}
The flowchart in Figure \ref{fig:flowchart} 
resumes the main steps of the proposed ROM approach.


\subsection{Lifting function method}\label{sec:lift}

Let us suppose to have only one inlet $\Gamma_i$, where nonhomogeneous boundary conditions for the velocity are imposed, and one outlet $\Gamma_o$, where instead we enforce nonhomogeneous boundary conditions for the pressure. We denote with  $\bm\chi^{\bm u}(\bm x)$ and $\bm\chi^{p}(\bm x)$ the lifting functions for velocity and pressure, respectively. 
Then each snapshot is modified as follows:
\begin{equation}
\bm u' (\bm x, t) = \bm u (\bm x, t) - 
g^{\bm u}(t) \bm\chi^{\bm u}(\bm x),
\label{shiftu}
\end{equation}
\begin{equation}
 p' (\bm x, t) =  p (\bm x, t) - 
 g^{p}(t) \bm\chi^{p}(\bm x),
 \label{schiftp}
\end{equation}
where $g^{\bm u}(t)$ is the prescribed boundary condition for the velocity on $\Gamma_i$  
and $g^{p}(t)$ is the pressure outflow on $\Gamma_{o}$ computed by \eqref{wind_discr}. 

In order to compute the lifting functions $\bm\chi^{\bm u}$ and $\bm\chi^{p}$ we solve a potential flow problem: 
\begin{equation}
    \begin{cases}
    \begin{aligned}
        \Delta \bm\chi^{\bm u} &= 0  \quad \mbox{ on } \Omega,\\ 
        \bm\chi^{\bm u} &= 1  \quad \mbox{ on } \Gamma_i,\\
        \bm\chi^{\bm u} &= 0  \quad \mbox{ on } \Gamma_w,\\
        \nabla \bm\chi^{\bm u} \cdot \bm n &= 1  \quad \mbox{ on } \Gamma_o,
    \end{aligned}
    \end{cases}\qquad \mbox{and} \qquad 
    \begin{cases}
    \begin{aligned}
        \Delta \bm\chi^p + \nabla\cdot ( \nabla \cdot (\bm\chi^{\bm u} \otimes \bm\chi^{\bm u})) &= 0   \quad \mbox{ on } \Omega, \\  
        \bm\chi^{p} &= 1  \quad \mbox{ on } \Gamma_o,\\
        \nabla \bm\chi^{p} \cdot \bm n &= 0  \quad \mbox{ on } \partial \Omega \setminus \Gamma_o.
    \end{aligned}
    \end{cases}
    \label{potentialflow}
\end{equation}

Then we apply the POD algorithm to the homogenized snapshots $\bm u'$ and $p'$, which satisfy homogeneous boundary conditions, to compute the reduced space, and the solution can be approximated as:
\begin{equation}
    \bm u \approx 
    g^{\bm u}(t) \bm\chi^{\bm u}(\bm x)+ \sum_{i=1}^{N_{\bm u}} a_i(t) \bm\phi_i(\bm x),
\end{equation}
\begin{equation}
    p \approx 
    g^{p}(t) \bm\chi^{p}(\bm x) + \sum_{i=1}^{N_p} b_i(t) \psi_i(\bm x).
\end{equation}
Note that $g^{p}(t)$ for time instances different from the training ones 
is the output of a neural network  properly employed to find a map between the training time instances and the corresponding pressure outflow values (see Sec. \ref{sec:nn}).

When a nonhomogeneous boundary condition for the pressure is imposed on more than one outlet section, for each outer boundary $j$, 
$\bm\chi^{p}_j$ is computed from the second system of \eqref{potentialflow} with 
$\bm\chi^p_j=1$ on the boundary $j$ and homogeneous Neumann conditions on all the other boundaries. In this case equation \eqref{schiftp} is generalized as follows: 

\begin{equation}
 p' (\bm x, t) =  p (\bm x, t) - 
\sum_{j=1}^{N^{p}_{BC}} g^{p}_j(t) \bm\chi^{p}_j(\bm x),
\end{equation}
where $N^{p}_{BC}$ denote the numbers of boundaries where nonhomogeneous Dirichlet boundary conditions are imposed on the pressure field. 

\subsection{Proper orthogonal decomposition}\label{sec:pod}
Let us consider the discrete finite dimensional set $\{ t_1, \dots, t_{{N_t}} \} \subset (t_0,T]$ of time instances, with cardinality $N_t$. The corresponding high fidelity solutions, for each field of interest $\Phi=\{p,\bm u\}$, are stored in the snapshot matrix $S_{\Phi}\in \mathbb{R}^{N_h\times N_t}$ as follows:

\begin{equation*}
     \mathcal{S}_{\Phi} = 
     \begin{Bmatrix}
    \Phi(t_1),  \dots, \Phi(t_{N_t})
    \end{Bmatrix}.
\end{equation*}
The POD space of size $N_{\Phi}$ consists of the set of basis $\mathcal{V} = \{\bm\zeta_1,\dots,\bm\zeta_{N_{\Phi}}\}$ that is the solution of the following minimization problem:
\begin{equation}
    \min_{\mathcal{V}} \Vert \mathcal{S}_{\Phi}- \mathcal{V}\mathcal{V}^T\mathcal{S}_{\Phi} \Vert \quad s.t. \quad \mathcal{V}^T\mathcal{V}=\mathcal{I},
    \label{min}
\end{equation}
i.e., the set of basis minimizing the distance between the snapshots and their projection onto the POD space. The POD procedure employes a least-square approach to squeeze data and find an optimal orthonormal basis \cite{bang2004greedy,hesthaven2018non}. 
Problem \eqref{min} is equivalent to the following eigenvalue problem \cite{kunisch2002galerkin}:
\begin{equation*}
    C_{\Phi} \bm c^{\Phi}_s = \lambda_s^{\Phi} \bm c_s^{\Phi}, \quad s = 1, \dots, N_t,
\end{equation*}
where $C_{\Phi}=\frac{1}{N_t}S_{\Phi}^TS_{\Phi}\in \mathbb{R}^{N_t\times N_t}$ is the correlation matrix associated with the snapshots. 

For a chosen $N_{\Phi} \ll \text{min}(N_{h},N_t)$, the reduced orthonormal basis is computed as follows:
\begin{equation*}
    \bm\phi_i=\frac{1}{\sqrt{\lambda^{\Phi}_i}} S_{\Phi}\bm c^{\Phi}_i, \quad i = 1,\dots, N_{\Phi}.
\end{equation*} 
A common approach to select the dimension of the reduced space $N_{\Phi}$ is to fix a threshold $\delta$  such that:
\begin{equation}
\frac{\sum_{i=1}^{N_{\Phi}}\lambda_i^{\Phi}}{\sum_{i=1}^{N_t}\lambda_i^{\Phi}} \ge \delta,
\label{eq:energy}
\end{equation}
where the left hand side is the percentage cumulative energy of the eigenvalues retained by the first $N_{\Phi}$ modes. The sum of the neglected singular values is a measure of the error introduced by ROM \cite{quarteroni2015reduced}:
\begin{equation}
\sum_{i=N_{\Phi}+1}^{\text{min}(N_{h},N_t)} \lambda_i^{\Phi}.
     \label{err}
\end{equation}
So, based on the value of $N_{\Phi}$, 
the ROM solution is enabled to approximate the FOM solution with arbitrary accuracy. 

\subsection{Galerkin projection}\label{sec:galerkin}

The $L^2$ orthogonal projection of equation \eqref{eq:NS1} 
onto  $\text{span}\{\bm\phi_1,\dots,\bm\phi_{N_{\bm u}}\}$ results in \cite{akhtar2009stability,bergmann2009enablers,lorenzi2016pod}:
\begin{equation}
    \big(\bm \phi_i, \partial_t \bm{u} + \nabla \cdot (\bm{u} \otimes \bm{u}) - \nabla \cdot (\nu \nabla \bm{u})+\nabla p \big)_{L^2(\Omega)} = 0 \quad \text{for }i = 1, \dots, N_{\bm u}.
    \label{proj_l2}
\end{equation}
By substituting equations \eqref{u_rb} and \eqref{p_rb} in equation \eqref{proj_l2} and given the orthonormality of the reduced basis, we get the following dynamical system: 
\begin{equation}
    \dot{\bm a} = \nu \bm B \bm a - \bm a^T \bm C \bm a - \bm K \bm b,
    \label{galerkin_u}
\end{equation}
where $\bm a = \{ a_i(t)\}_{i=1}^{N_u}$ and $\bm b = \{ b_i(t)\}_{i=1}^{N_p}$ and 
\begin{equation}
    B_{ij} = \big( \bm\phi_i, \Delta \bm\phi_j \big)_{L^2(\Omega)},
\end{equation}
\begin{equation}
    C_{ijk} = \big( \bm\phi_i, \nabla \cdot (\bm\phi_j \otimes \bm\phi_k) \big)_{L^2(\Omega)},
\end{equation}
\begin{equation}
    K_{ij} = \big( \bm\phi_i, \nabla \psi_j \big)_{L^2(\Omega)}.
\end{equation}\\[2ex]
The continuity equation \eqref{eq:NS2} is projected as well onto $\text{span}\{\psi_1, \dots, \psi_{N_p}\}$ by providing: 
\begin{equation}
    \big( \psi_i,\nabla \cdot \bm u \big)_{L^2(\Omega)} = 0 \quad \text{for }i = 1, \dots, N_{p}.
    \label{proj_p_continuity}
\end{equation}
By substituting  equation \eqref{u_rb} in equation \eqref{proj_p_continuity}, we get:
\begin{equation}
    \bm P {\bm a} = 0,
    \label{galerkin_p_continuity}
\end{equation}
where
\begin{equation}
    P_{ij} = \big( \bm\psi_i, \nabla \cdot\bm\phi_j \big)_{L^2(\Omega)}.
\end{equation}
It is well known that the numerical approximation of incompressible Navier-Stokes equations shows stability issues. In particular, the discrete function spaces for velocity and pressure need to satisfy the inf-sup condition \cite{boffi2013mixed,brezzi1990discourse}. At full order level, the finite element method requires proper basis functions for pressure and velocity (such as the Taylor-Hood elements) \cite{schulz2009one,gunzburger2002perspectives,ali2020stabilized}. However, the finite volume method does not require any particular treatment because no projection of the equations is required.

On the other hand, at reduced order level, we perform a Galerkin project requiring to ensure the reduced version of the inf-sup condition \cite{stabile2018finite, ballarin2015supremizer,gerner2012certified,rozza2007stability}. At this aim, in this work we consider two different strategies: the supremizer enrichment and the Pressure Poisson equation.
\\[2ex]
\emph{\bf Supremizer enrichment: }
A popular method provided in literature to employ the reduced stability condition is the supremizer enrichment of the velocity space. For each pressure basis function $\psi_i$, a so-called supremizer $\bm s_i$ is computed and added to the velocity reduced space:
\begin{equation}
    \mathcal{V} = \{\bm\phi_1,\dots,\bm\phi_{N_{\bm u}},\bm s_1,\dots,\bm s_{N_p} \}.
\end{equation}
For each  $\psi_i$ the corresponding supremizer element $\bm s_i$ can be computed by solving the following system: 
\begin{equation}
    \begin{cases} 
        \Delta \bm s_i =  -\nabla \psi_i & \quad \mbox{ on }  \Omega, \\ 
        \bm s_i = 0 & \quad \mbox{ on }  \partial \Omega .
    \end{cases} 
    \label{supremizer}
\end{equation}
This is the so-called exact supremizer enrichment. However, also the so-called approximate supremizer enrichment, where the pressure basis in \eqref{supremizer} is replaced by the pressure snapshots, has shown promising results: see, e.g., \cite{stabile2018finite}. 
Further details about the  supremizer enrichment method, both in its exact and approximate formulation, can be found in \cite{ballarin2015supremizer,gerner2012certified,rozza2007stability}. \\[2ex]
\emph{\bf Pressure Poisson equation}: Another popular approach to solve the stability issues of the  reduced problem is based on the employment of a Poisson equation for the pressure in place of the continuity equation \eqref{eq:NS2}. The Poisson equation for the pressure can be obtained by taking the divergence of the momentum equation \eqref{eq:NS1} and by exploiting the divergence-free constraint \eqref{eq:NS2}. Neumann boundary conditions for the pressure are obtained by projecting the momentum equation \eqref{eq:NS1} onto the surface normal vector $\bm n$. Further details can be found in \cite{lorenzi2016pod}.

At the full order level, one obtains the following system: 
\begin{equation}
    \begin{cases} 
    \Delta p = - \nabla \cdot ( \nabla \cdot (\bm{u} \otimes \bm{u})) & \quad \mbox{ on } \Omega, \\ 
    \frac{\partial p}{\partial \bm n} = -\nu \bm n \cdot (\nabla \times \nabla \times \bm u)   & \quad \mbox{ on } \partial \Omega .
     \end{cases}\label{poission} 
\end{equation}
Then equation \eqref{poission} is projected onto $\text{span}\{\psi_1, \dots, \psi_{N_p}\}$:
\begin{equation}
    \big(\nabla \psi_i,\nabla p \big)_{L^2(\Omega)} - \big(\psi_i,\nabla p \cdot \bm n \big)_{L^2(\Gamma_o)}= \big( \psi_i, \nabla \cdot (\nabla \cdot (\bm{u} \otimes \bm{u})) \big)_{L^2(\Omega)}.
    \label{proj_p}
\end{equation}
Let us replace $\bm u$ and $p$ with $\bm u_{\text{rb}}$ and $p_{\text{rb}}$ (see equations \eqref{u_rb} and \eqref{p_rb}). Then equation \eqref{proj_p} can be reformulated as follows:
\begin{equation}
    \bm D \bm b - \bm N \bm b = \bm a^T \bm G \bm a,
    \label{galerkin_p}
\end{equation}
where 
\begin{equation}
    D_{ij}=\big( \nabla \psi_i, \nabla \psi_j \big)_{L^2(\Omega)},
\end{equation}
\begin{equation}
    N_{ij}=\big( \psi_i, \nabla \psi_j \bm n \big)_{L^2(\Gamma_o)},
\end{equation}
\begin{equation}
    G_{ijk}=\big( \psi_i, \nabla \cdot(\nabla \cdot (\bm \phi_j\otimes\bm\phi_k)) \big)_{L^2(\Omega)}.
\end{equation} \\
So the reduced system is given by equations \eqref{galerkin_u} and \eqref{galerkin_p}.
\textcolor{black}{The computational burden associated with nonlinear terms, $\bm C$ and $\bm G$, can be addressed by keeping under control the number of basis functions. In fact, notice that the size of tensor $\bm C$ scales as $N_{\bm u}^3$, while the size of tensor $\bm G$ scales as $N_{\bm u}^2 N_p$.}

The choice between supremizer enrichment and PPE to stabilize ROMs is problem-dependent. In \cite{stabile2018finite}, they found that the PPE approach is more effective for long-time integration. However, experimental validation on a specific problem ultimately identifies the most appropriate method. This topic is explored further in Section~\ref{sec:cylinder}.

\subsection{Neural network interpolation}\label{sec:nn}

In the offline phase, we train a feed-forward neural network for approximating the temporal dynamics of the pressure outlet boundary conditions: 
\begin{equation}\label{eq:pi}
    g^p_b: t_i \mapsto \left[g^p_b(t_i)\right]_{b=1}^{N^p_{BC}}, \quad i=1,\dots,N_T.
\end{equation}
Once the network is trained, the outflow pressure can be computed online for every new time value $t_{\text{new}}$ in the online phase which is not originally included in the full order database, i.e. for every test point.

\begin{figure}[htb]
\centering
\begin{tikzpicture}[x=2.4cm,y=1.2cm,scale=0.77]
  \readlist\Nnod{2,4,4,2} 
  \readlist\Nstr{n,m,k} 
  \readlist\Cstr{i,h^{(\prev)},o} 
  \def\yshift{0.55} 
  
  \foreachitem \N \in \Nnod{
    \def\lay{\Ncnt} 
    \pgfmathsetmacro\prev{int(\Ncnt-1)} 
    \foreach \i [evaluate={\c=int(\i==\N); \y=\N/2-\i-\c*\yshift;
                 \x=\lay; \n=\nstyle;
                 \index=(\i<\N?int(\i):"\Nstr[\n]");}] in {1,...,\N}{ 
      \node[node \n] (N\lay-\i) at (\x,\y) {$\strut\Cstr[\n]_{\index}$};
      
      \ifnumcomp{\lay}{>}{1}{ 
        \foreach \j in {1,...,\Nnod[\prev]}{ 
          \draw[white,line width=1.2,shorten >=1] (N\prev-\j) -- (N\lay-\i);
          \draw[connect] (N\prev-\j) -- (N\lay-\i);
        }
        \ifnum \lay=\Nnodlen
          \draw[connect] (N\lay-\i) --++ (0.5,0); 
        \fi
      }{
        \draw[connect] (0.5,\y) -- (N\lay-\i); 
      }
      
    }
    \path (N\lay-\N) --++ (0,1+\yshift) node[midway,scale=1.6] {$\vdots$}; 
  }
  
  \node[above=3,align=center,mydarkgreen] at (1,-3) {Input\\[-0.2em]layer};
  \node[above=2,align=center,mydarkblue] at (2.5,-1) {Hidden\\[-0.2em]layers};
  \node[above=3,align=center,mydarkred] at (4,-3) {Output\\[-0.2em]layer};
  
\end{tikzpicture}

\flushleft\caption{Sketch of a feedforward neural network.}
\label{my_net}
\end{figure}
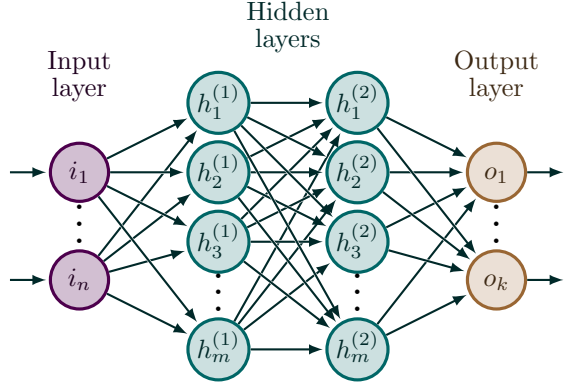

For sake of completeness, we are going to briefly describe how feedforward neural networks work. For further details, the reader is referred, e.g., to \cite{goodfellow2016deep,kriesel2007brief,calin2020deep}. A fully connected feedforward neural network is an architecture composed of a set of neurons (or nodes) arranged in layers, and each neuron in a certain layer is connected to all the neurons in the next layer through oriented edges (or synapses) \cite{rosenblatt1958perceptron,minsky1969introduction,fine2006feedforward}.
The number of neurons in the input and output layers is related to the specific problem the network solves.  Other neurons make up the hidden layers and their number is chosen to improve the results (see Figure \ref{my_net}).

Each neuron $j$ is identified by three functions:
\begin{itemize}
    \item the propagation function $u_j$: 
    $$u_j=\sum_{k=1}^m w_{s_k,j}y_{s_k} + \beta_j,$$ where $\beta_j$ is the bias, $y_{s_k}$ is the output related to the sending neuron $k$,  $w_{s_k,j}$ are the weights and $m$ is the number of sending neurons connected with the neuron $j$.
    \item the activation function $a_j$: 
    $$ a_j=f_{\text{act}}\left(\sum_{k=1}^m w_{s_k,j}y_{s_k}+\beta_j\right).$$
    Possible choices are the sigmoid function, hyperbolic tangent, RELU, SoftMax, or the Softplus activation function \cite{sharma2017activation}, the preferred choice in our application.
    \item the output function $y_j$:
    \begin{equation*}
        y_j=f_{\text{out}}(a_j). 
    \end{equation*}
    Often the output function coincides with the identity, so that $a_j = y_j$.
\end{itemize}

The backpropagation algorithm is engaged to optimize the weights of the synapses \cite{rumelhart1986learning,rojas1996backpropagation}. A loss function $\mathcal{L}$, which quantifies the distance between the actual output and the desired output, is minimized by evaluating (backward) the gradient relative to the weights. The model parameters are tuned as follows:
\begin{equation}
\begin{split}
    &  \bm{w}=\bm{w}-\eta \frac{\partial \mathcal{L}}{\partial \bm{w}}, \\
    & \bm{b} = \bm{b} -\eta  \frac{\partial \mathcal{L}}{\partial \bm{b}},
\end{split}
\end{equation}
where $\bm{w}$ and $\bm b$ are matrix representations of the weights and biases and $\eta$ is the learning rate. 

Hyperparameters like the activation function, the number of layers, the number of neurons per layer, and the learning rate are adapted to enhance the learning process. For more details, the reader is referred to \cite{sharma2017activation,montana1989training}. 



\section{Numerical results}\label{sec:results}
In order to validate our ROM approach, we consider two different test cases:
\begin{itemize}
    \item \textit{Case 1}: an idealized blood vessel consisting of a 2D channel flow. Despite its semplicity, this academic benchmark allows us to highlight the basic performance of our ROM approach. 
    \item \textit{Case 2}: a patient-specific thoracic aortic arch, equipped with five outlet sections. This example highlights the capability of our ROM framework to handle complex cases. 
\end{itemize}

All the FOM simulations are run in \openfoam (\url{https://openfoam.org/}) 
, whereas the home-made library ITHACA-FV (\url{https://github.com/ITHACA-FV/ITHACA-FV}) is employed for the ROM computations \cite{stabile2018finite}. PyTorch (\url{https://pytorch.org/}) is the library employed for the training of the neural network.

\subsection{Case 1: Idealized blood vessel}
\label{sec:cylinder}
\begin{figure}
    \centering
    \includegraphics[width=.35\textwidth]{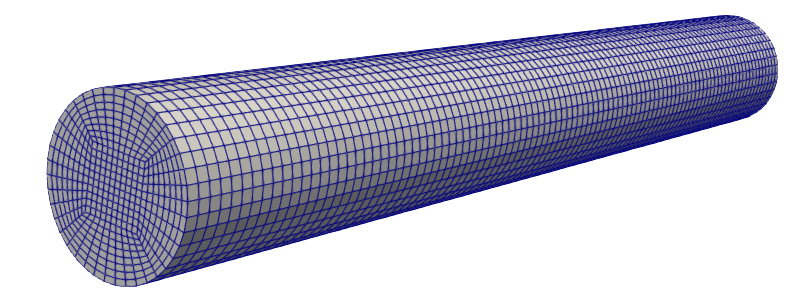}
    \caption{\textit{Case 1}: sketch of the domain and of the computational grid.}
    \label{fig:domain_mesh}
\end{figure}

The dynamics of the blood flow in a cylinder over the time window $(t_0, T] = (0, 1]$ s is here investigated. 
The radius and length are $R=2$ cm and $L=30$ cm, respectively. 
The mesh has been generated by using the \emph{blockmesh} utility available in \openfoam. The geometry as well as the mesh are shown in Figure \ref{fig:domain_mesh}. The number of cells, the minimum and the maximum size of the mesh, the maximum non-orthogonality and the maximum skewness are reported in Table \ref{tab:features_mesh}. 

\begin{figure}
	\centering
    \subfloat[][Inflow velocity.\label{fig:BCU}]{\includegraphics[width=.4\textwidth]{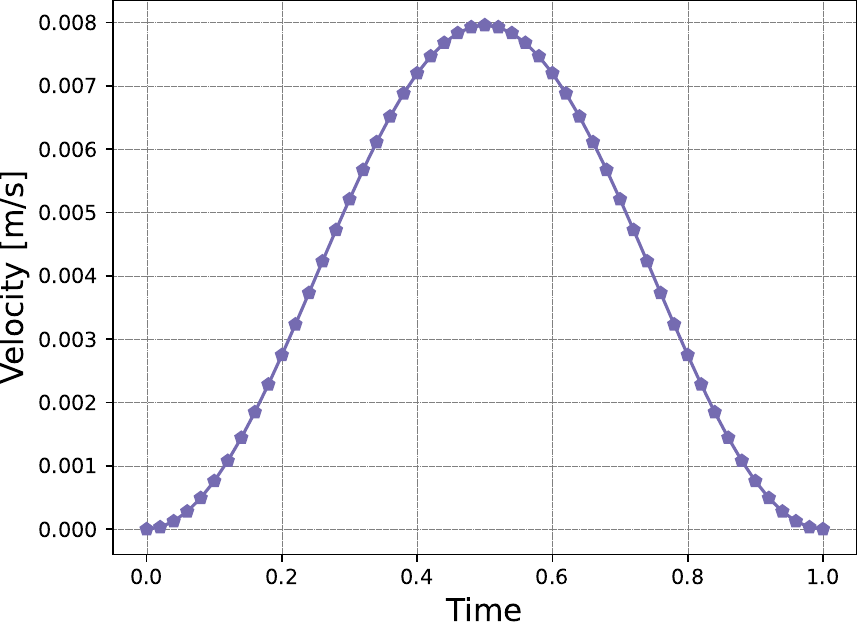}}
    \hspace{2ex}
 	\subfloat[][Outflow pressure.\label{fig:BCp}]{\includegraphics[width=.4\textwidth]{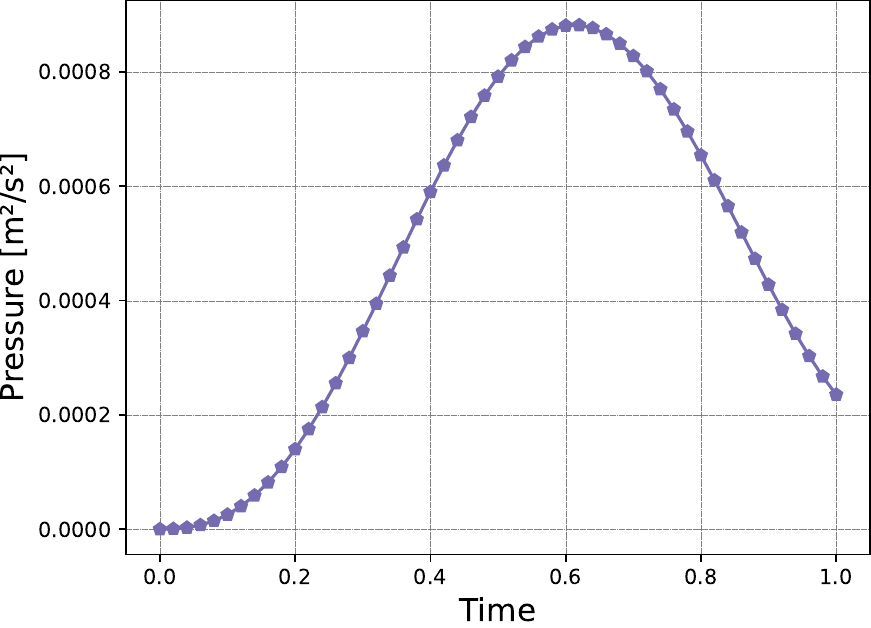}}\\
	\caption{\textit{Case 1}: Time evolution of the boundary conditions for the velocity (a) and the pressure (b).}
	\label{BC}
\end{figure}

\begin{table}[h]
\caption{\textit{Case 1}: Features of the mesh.}
\centering
\renewcommand{\arraystretch}{1.5}
\begin{tabular}{ccccc}
\hline
\rowcolor{gray!20} 
 Number of cells & $h_{\text{min}}$ (m) & $h_{\text{max}}$ (m) & max non-orthogonality ($^\circ$) & max skewness\\
\hline
42000 & 1.68$\cdot 10^{-3}$ & 2.42$\cdot 10^{-3}$ & 29.9 & 0.8612 \\
\hline
\end{tabular}
\label{tab:features_mesh}
\end{table}
The unsteady profile  with amplitude $u_0$ 
\begin{equation}\label{eq:sinprofile}
    u(t) = u_0 \sin^2\left(\frac{t}{4R_dC}\right)
\end{equation}
is imposed on the inlet, whereas the outflow pressure is governed by the Windkessel model \eqref{eq:Windkessel-FOM}. In this case, an analytical solution can be obtained:
\begin{equation}\label{eq:pBC}
    p(t) = {R_d  u_0 \pi R^2}\left[ \left(\frac{R_p}{R_d}+\frac{1}{2}\right)\sin^2\left(\frac{t}{4R_dC}\right)  + \frac{1}{4} \left(1-e^{\frac{t}{2R_dC}} - \sin\left(\frac{t}{2R_dC}\right)\right)\right].
\end{equation}
The kinematic viscosity of the fluid is $\nu=0.004$ m$^2$/s. We set $u_0 = 0.007957$. 
The Windkessel coefficients are taken from \cite{moghadam2013modular}. For sake of completeness they are reported in Table \ref{tab:parameters}.
The time dependent boundary values for velocity \eqref{eq:sinprofile} and pressure \eqref{eq:pBC} are shown in Figure \ref{BC}. 
\begin{table}[h]
\caption{\textit{Case 1}: Parameter values in equations \eqref{eq:sinprofile}  and \eqref{eq:pBC}.} 
\centering
\renewcommand{\arraystretch}{1.5}
\begin{tabular}{ccc}
\hline
\rowcolor{gray!20} 
  $R_p$ [m$^{-1}$s$^{-1}$] & $R_d$ [m$^{-1}$s$^{-1}$] & C [ms$^{2}$] \\
\hline
 $10^4$ & $10^5$ & 0.07957$\cdot 10^{-5}$  \\
\hline
\end{tabular}
\label{tab:parameters}
\end{table}

The FOM simulation is run with a time step $\Delta t = 5\cdot10^{-5}$ s. This value has been selected to ensure the stability of the simulation; moreover, we have verified that the solution exhibits radial symmetry. To train our ROM, we collect a total of 50 time-dependent snapshots over the time interval (0,1] s, one every 400 time step (i.e. one every 0.02 s).  For this test case, we do not consider any new time instance in the online phase, i.e. we run our ROM over the time interval (0,1] s with a reduced time step $\Delta t_r=0.02$ s. Consequently, no neural network for mapping the outflow pressure is introduced. 
In the following, unless specified otherwise, the supremizer approach is adopted to stabilize the ROM. In addition, the same number of modes $N_{\Phi}$ is employed for velocity, pressure, and supremizers.

\begin{figure}
	\centering
    \subfloat[][\label{fig:error_Nu}]{\includegraphics[width=.4\textwidth]{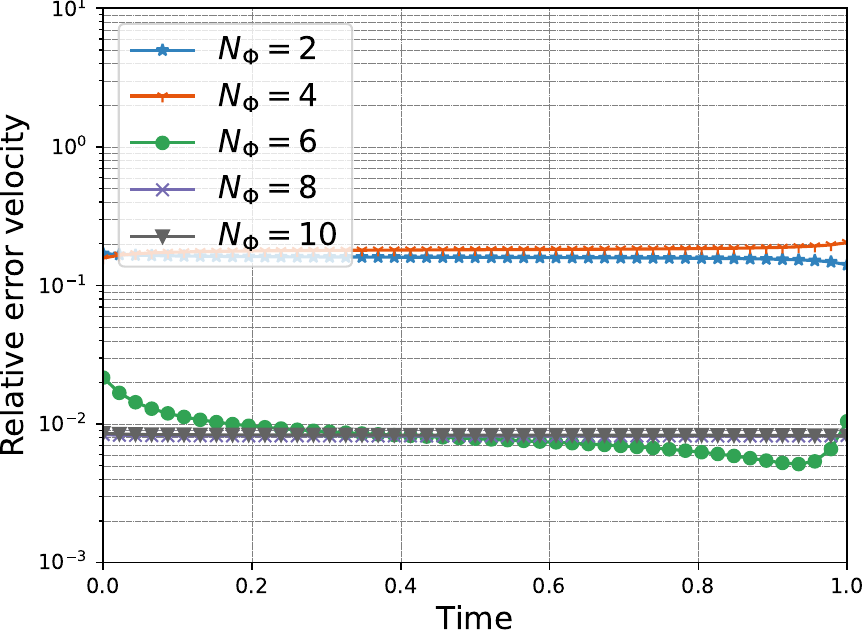}}\hspace{2ex}
 	\subfloat[][\label{fig:error_Np}]{\includegraphics[width=.4\textwidth]{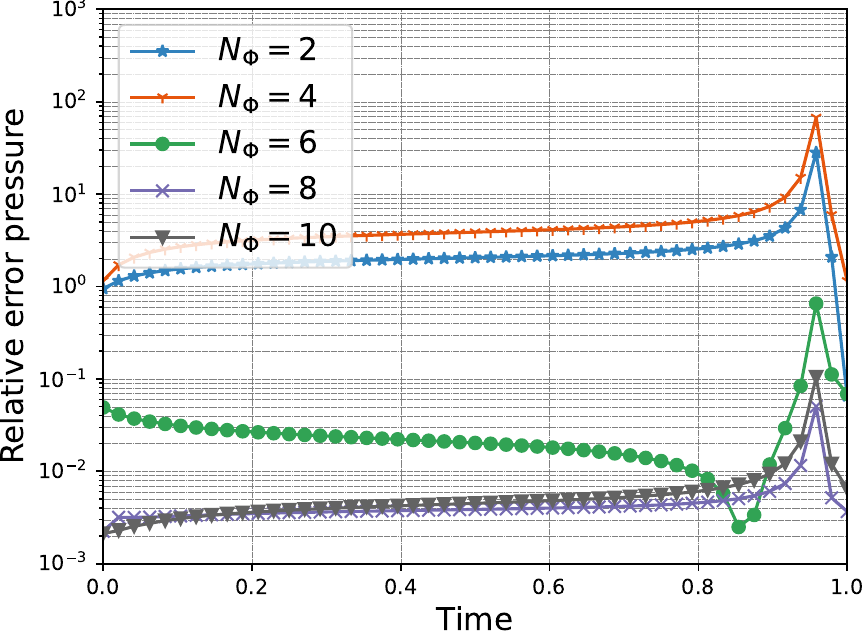}}\\
    
    \subfloat[][\label{fig:error_relative_U_PPE_N}]{\includegraphics[width=.4\textwidth]{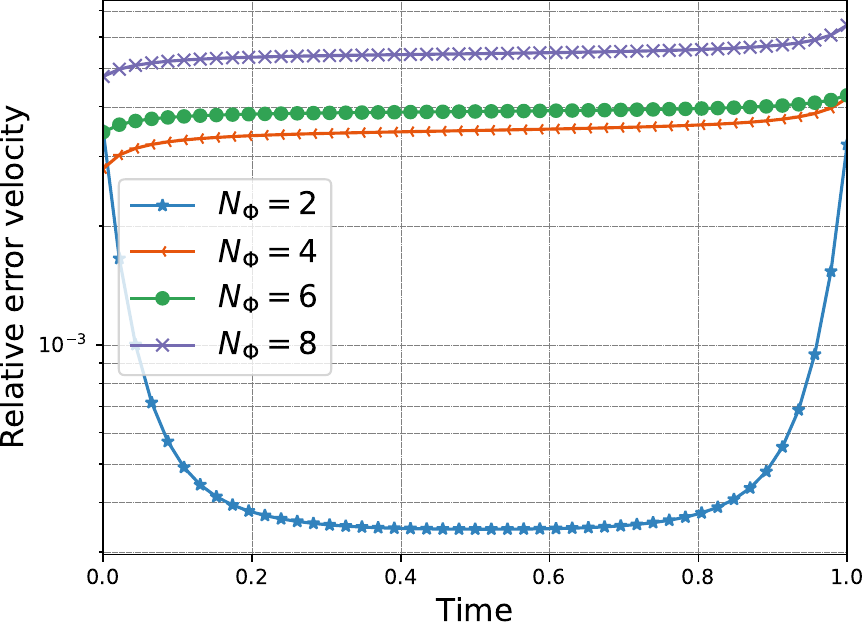}}\hspace{2ex}
    \subfloat[][\label{fig:error_relative_P_PPE_N}]{\includegraphics[width=.4\textwidth]{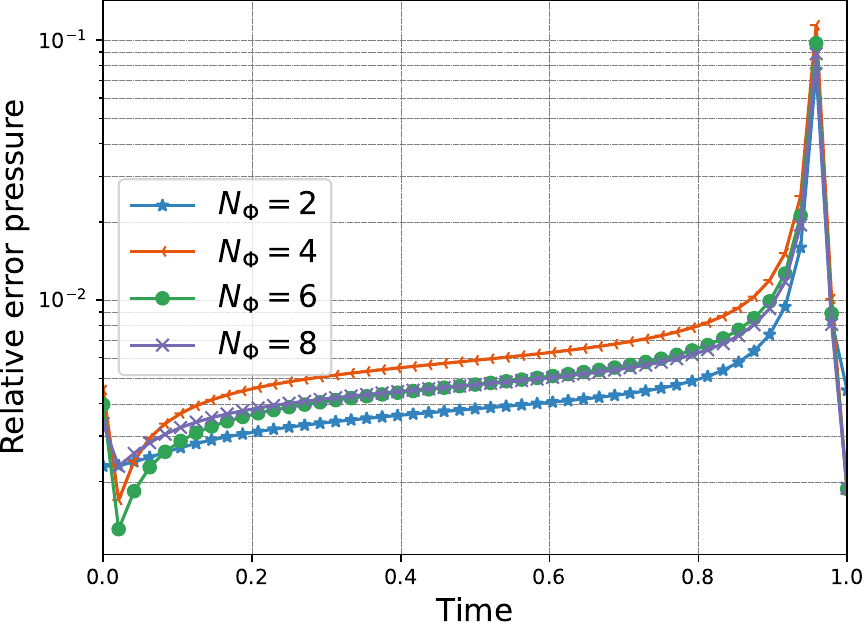}}
	\caption{\textit{Case 1}: Time evolution of the relative reconstruction error for velocity and pressure as the number of modes $N_{\Phi}$ increases. The supremizer approach is adopted in (a) and (b) while the PPE approach is adopted in (c) and (d). The number of modes $N_{\Phi}$ is the same for pressure, velocity and supremizers.}
	\label{error_N}
\end{figure}

\begin{figure}
	\centering
    \subfloat[][\label{fig:error_proju}]{\includegraphics[width=.4\textwidth]{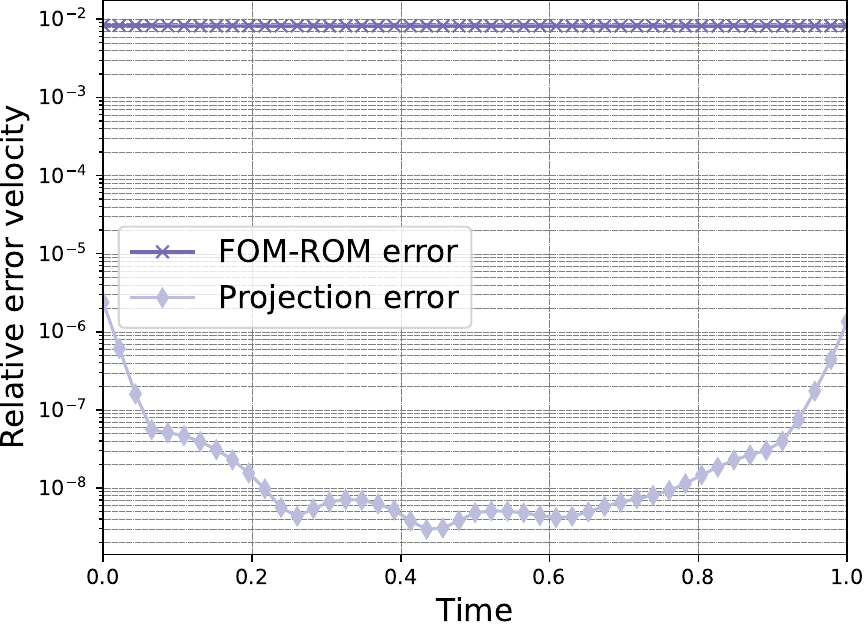}}\hspace{2ex}
 	\subfloat[][\label{fig:error_projp}]{\includegraphics[width=.4\textwidth]{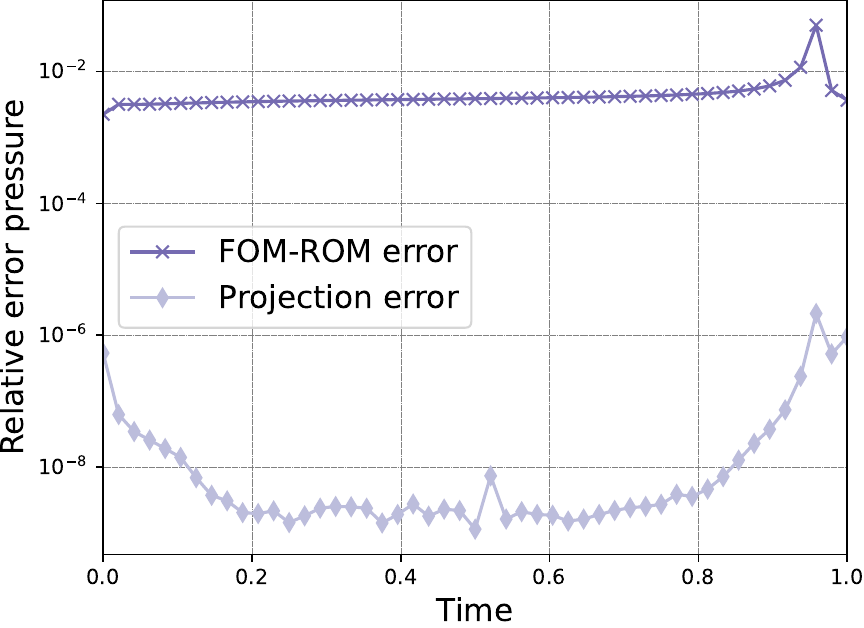}}\\
	\caption{\textit{Case 1}: Comparison between reconstruction error and projection error for velocity and pressure with $N_{\Phi}=8$. The supremizer approach is adoped and $N_{\Phi}$ is the same for pressure, velocity and supremizers.}
	\label{error_proj}
\end{figure}

\begin{figure}[!htb]
    \centering
    \begin{minipage}{.7\textwidth}
        \centering
        \begin{overpic}[width=0.8\textwidth]{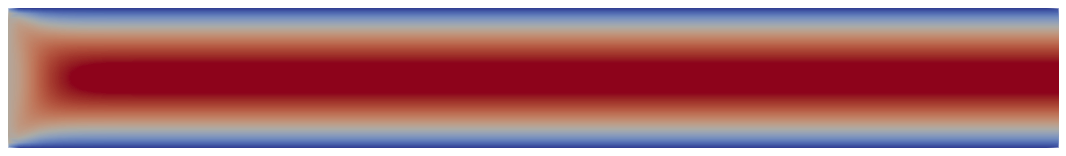} 
        \put(45,15){$\bm u$ FOM}
        \end{overpic}\\ \vspace{3ex}
        \begin{overpic}[width=0.8\textwidth]{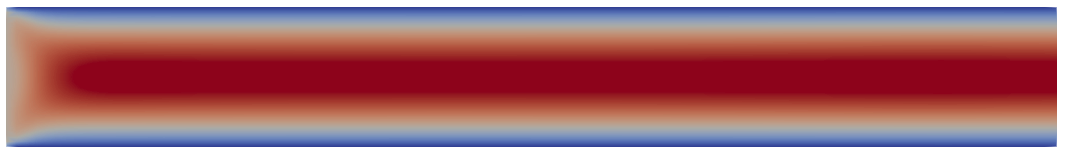} 
        \put(45,15){$\bm u$ ROM}
        \end{overpic}
    \end{minipage}%
    \begin{minipage}{0.29\textwidth}
        \centering
        \includegraphics[width=0.4\linewidth]{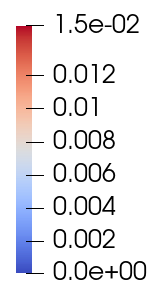}
    \end{minipage}
\caption{\textit{Case 1}: Qualitative comparison of FOM-ROM velocity at $t=0.5$ s and $N_\phi=8$.}
\label{fig:u-fom-rom}
\end{figure}

We have that 2 POD modes are enough to get over $99.99\%$ of the cumulative energy of the eigenvalues (see equation \eqref{eq:energy}) for both pressure and velocity. The time evolution of the relative reconstruction errors between FOM and ROM for the velocity and pressure fields, 
\begin{equation}
    \varepsilon_u = \frac{\| \bm  u_{\text{FOM}}- \bm u_{\text{ROM}}\|_{L^2(\Omega)}}{\|\bm u_{\text{FOM}}\|_{L^2(\Omega)}} \qquad \text{and} \qquad \varepsilon_p = \frac{\| p_{\text{FOM}}-  p_{\text{ROM}}\|_{L^2(\Omega)}}{\|p_{\text{FOM}}\|_{L^2(\Omega)}},
    \label{eq:errors-rel}
\end{equation}
is depicted in Figure \ref{fig:error_Nu} and \ref{fig:error_Np} at varying of the number of POD modes. 
In particular, we note that for $N_{\Phi}=6$  a relative error below 2\% for the velocity is obtained over the entire time interval. On the other hand, for the pressure error exhibits a more irregular trend: from $t = 0$ to $t = 0.94$ s  it remains below 10\% reaching a minimum value of $0.25\%$, but it becomes extremely large (about 65\%) towards the end of the simulation ($t=0.96$ s). We need to set $N_{\Phi}=8$ to significantly get down the maximum error for the pressure achieving the $3\%$. For $N_{\Phi} > 8$ the error does not exhibit significant changes. 
A comparison between the reconstruction error and the projection error (i.e., the error obtained by projecting the snapshots onto the basis chosen representing the best achievable result) is shown in Figure \ref{error_proj} for $N_\Phi = 8$. As expected, the projection error is smaller than the reconstruction error over the entire time window both for pressure and velocity. 
It is worth to note the presence of a small amplitude peak appears for the pressure at at the end of the time window, both in the reconstruction error and in the projection one (see Figure~\ref{fig:error_projp}). 
We also show a FOM-ROM qualitative comparison in Figure \ref{fig:u-fom-rom} for $t=0.5$ s. We see that our ROM is able to capture the main dynamics of the flow on the whole domain both for pressure and velocity. 

\begin{figure}[!htb]
    \centering
    \begin{minipage}{.7\textwidth}
        \centering
        \begin{overpic}[width=0.8\textwidth]{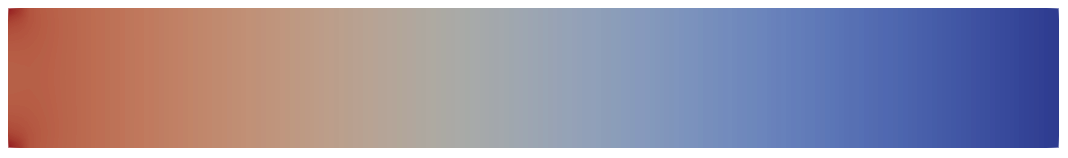} 
        \put(45,15){$p$ FOM}
        \end{overpic}\\ \vspace{3ex}
        \begin{overpic}[width=0.8\textwidth]{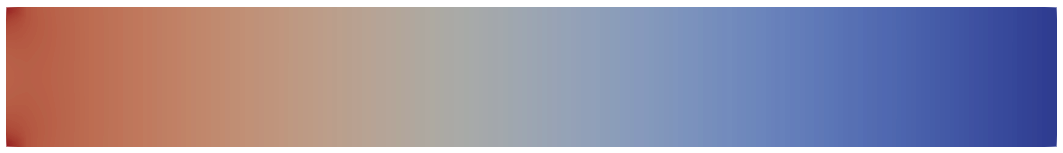} 
        \put(45,15){$p$ ROM}
        \end{overpic}
    \end{minipage}%
    \hfill
    \begin{minipage}{0.29\textwidth}
        \centering
        \includegraphics[width=0.4\linewidth]{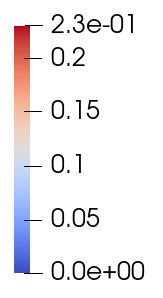}
    \end{minipage}
\caption{\textit{Case 1}: Qualitative comparison of FOM-ROM pressure at $t=0.5$ s and $N_\phi=8$.}
\label{fig:p-fom-rom}
\end{figure}

\begin{figure}
	\centering
    \subfloat[][ $N_{\Phi}=8$.\label{fig:error_relative_U_N_8_PPE_vs_sup}]{\includegraphics[width=.4111\textwidth]{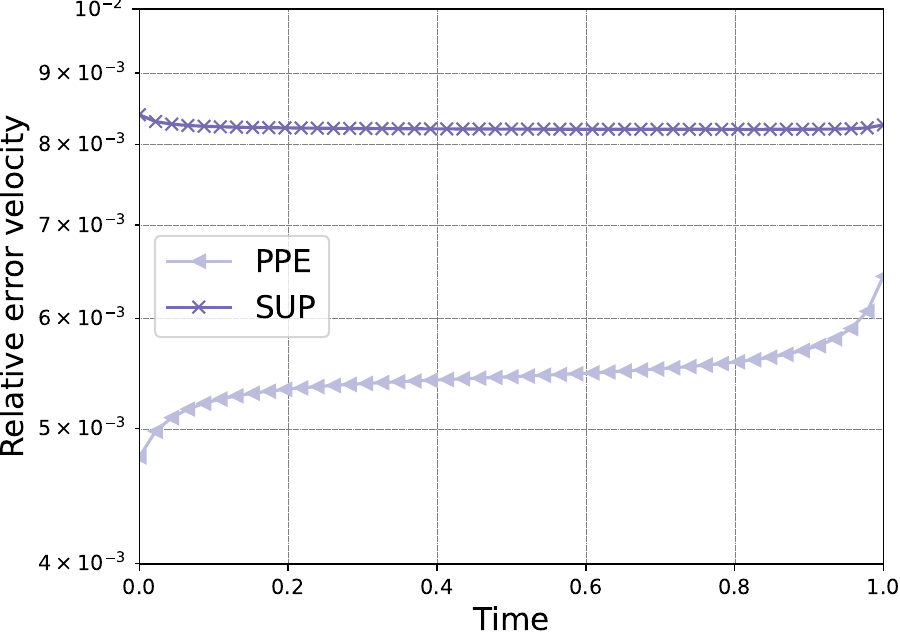}}\hspace{2ex}
 	\subfloat[][ $N_{\Phi}=8$.\label{fig:error_relative_P_N_8_PPE_vs_sup}]{\includegraphics[width=.4\textwidth]{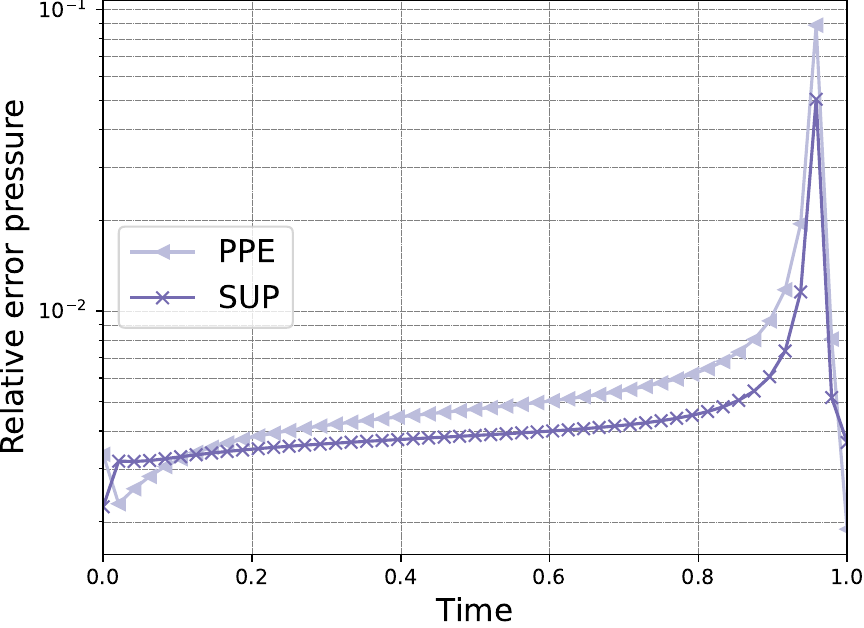}}

	\caption{\textit{Case 1}: Comparison in terms of the relative reconstruction error between PPE and supremizer (SUP) 
   for  $N_{\Phi}=8$.}
	\label{error_relative_PPE_N}
\end{figure}

Finally we consider also the PPE approach and  perform a comparison against the supremizer approach. 
In Figure \ref{fig:error_relative_U_PPE_N} and \ref{fig:error_relative_P_PPE_N}, the relative reconstruction error for pressure and velocity is computed as $N_{\Phi}$ increases and the PPE stabilization is adopted. Unlike the supremizer approach, we observe that the PPE error is not significantly affected by the number of the basis chosen. The pressure error remains around $1\%$ over the most part of the time interval for all the $N_{\Phi}$ values with a maximum value of 10\% just before the end of the simulation. Surprisingly, the velocity error increases as $N_{\Phi}$ rises but it remains bound in the order of $0.1\%$. Figure \ref{fig:error_relative_U_N_8_PPE_vs_sup} and \ref{fig:error_relative_P_N_8_PPE_vs_sup} show the velocity and pressure relative reconstruction error, respectively, for PPE and supremizer approach, with $N_{\Phi}=8$. We note that the velocity is better predicted with the PPE approach, 
although the difference is not so remarkable. On the contrary, the pressure is slightly better approximated with the supremizer approach almost everywhere. 
In correspondence of the peak before the end of the simulation the difference between the two approaches is of about $5\%$. 


\subsection{Case 2: Aortic arch}
Now we consider the dynamics of the blood flow in a patient-specific aortic arch.

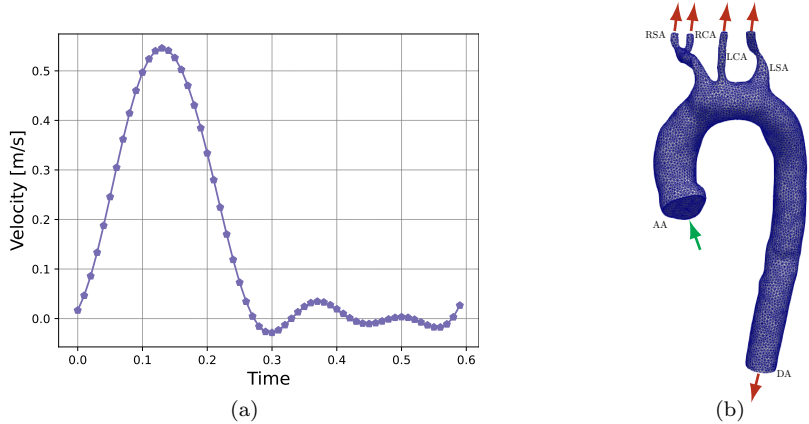
\begin{figure}
	\centering
    \subfloat[][\label{fig:BC_aorta}]{\includegraphics[width=.4\textwidth]{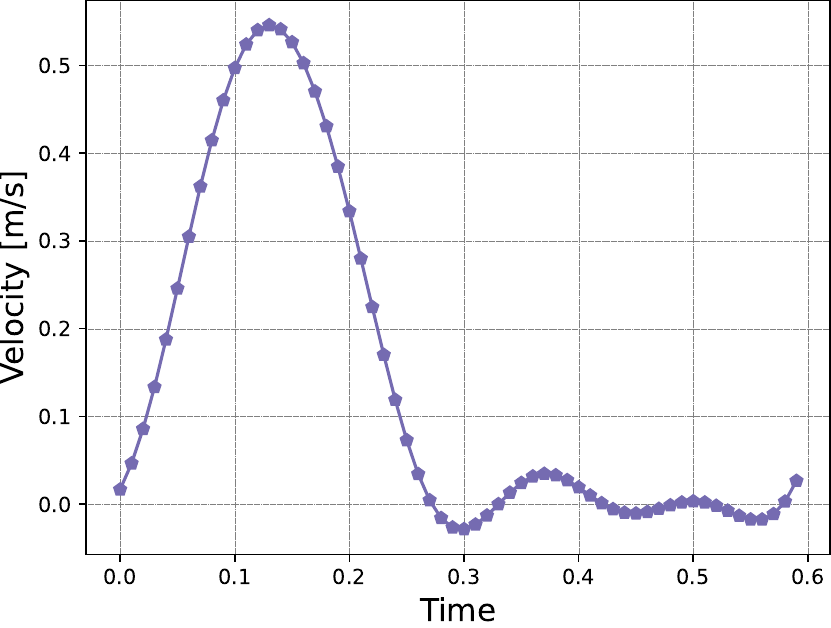}}\hspace{2cm}
 	\subfloat[][\label{fig:aorta_mesh}]{
    \begin{overpic}[width=0.15\textwidth]{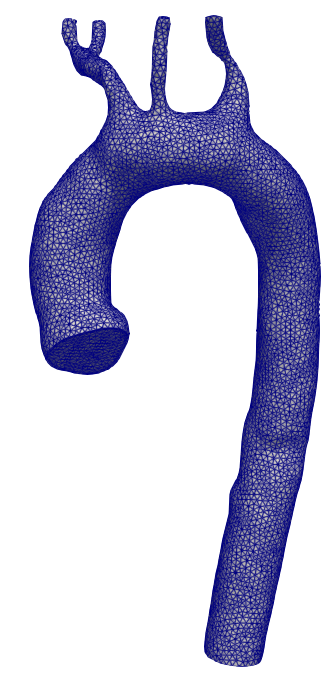}

    \put(4,44){\scalebox{0.38}{AA}}
    \put(38,3){\scalebox{0.38}{DA}}
    \put(2,96){\scalebox{0.38}{RSA}}
    \put(15.5,96){\scalebox{0.38}{RCA}}
    \put(24.2,90){\scalebox{0.38}{LCA}}
    \put(36,87){\scalebox{0.38}{LSA}}
    
    \linethickness{1pt}\put(17,38){\color{Green}\vector(-0.08,0.22){3}}
    
    \linethickness{1pt}\put(33,4){\color{BrickRed}\vector(-0.05,-0.20){2}}

     \linethickness{1pt}\put(14.5,98){\color{BrickRed}\vector(0.01,0.08){1}}

     \linethickness{1pt}\put(10,98){\color{BrickRed}\vector(0.01,0.08){1}}

     \linethickness{1pt}\put(23.5,98){\color{BrickRed}\vector(0.01,0.08){1}}

     \linethickness{1pt}\put(31,98){\color{BrickRed}\vector(0.01,0.08){1}}
    \end{overpic}
    }\\
	\caption{\textit{Case 2}: Time evolution of the boundary condition for the velocity (a) and sketch of the computational domain (b).}
	\label{fig:BC_aorta_mesh}
\end{figure}

The geometry in Figure \ref{fig:aorta_mesh}. The boundary includes an inlet section indicated with a green arrow,  Ascending Aorta (AA), and five outlet sections, Right Subclavian Artery (RSA), Right Common Carotid Artery (RCA), Left Common Carotid Artery (LCCA), Left Subclavian Artery (LSA) and Descending Aorta (DA),  indicated with red arrows.

The grid is generated with the open source mesh generator \emph{gmsh} (\url{https://gmsh.info/}). A mesh convergence analysis is performed in \cite{Girfoglio2020} for the same geometry. Therefore, for the current investigation, we adopt the finest mesh employed in \cite{Girfoglio2020}: see Figure \ref{fig:aorta_mesh}. 
The features of such a grid are reported in Table \ref{table_mesh}. 
\begin{table}[htb!]
\caption{\textit{Case 2}: Features of the mesh.}
\centering
\renewcommand{\arraystretch}{1.5}
\begin{tabular}{ccccc}
\hline
\rowcolor{gray!20} 
 Number of cells & $h_{\text{min}}$ (m) & $h_{\text{max}}$ (m) & mean non-orthogonality ($^\circ$) & max skewness\\
\hline
228296 & 5.8 $\cdot 10^{-4}$ & 3$\cdot 10^{-3}$ & 29.5 & 1.13 \\
\hline
\end{tabular}
\label{table_mesh}
\end{table}

A realistic time-dependent waveform is enforced on the ascending aorta for the velocity. It is shown in Figure~\ref{fig:BC_aorta} for a cardiac cycle of 0.6 s. 
The values of the Windkessel coefficients for each outlet are shown in Table \ref{table_windkessel_aorta}. Both the inflow boundary condition and the Windkessel coefficients are extracted from RHC tests and ECHO tests reported in \cite{Girfoglio2020}. The kinematic viscosity is $\nu=3.7\cdot 10^{-6}$ m$^2$/s.
\begin{table}[htb!]
\caption{\textit{Case 2}: Parameter values of the Windkessel model.} 
\centering
\renewcommand{\arraystretch}{1.5}
\begin{tabular}{cccc}
\hline
\rowcolor{gray!20} 
 Outlet & $R_p$ [m$^{-1}$s$^{-1}$] & $R_d$ [m$^{-1}$s$^{-1}$] & $C$ [ms$^{2}$] \\
\hline
Right subclavian artery (RSA) & 1.84$\cdot 10^{8}$  & 3.11$\cdot 10^{9}$ & 3.26$\cdot 10^{-5}$ \\
\hline
Right common carotid artery (RCA) & 1.23$\cdot 10^{8}$ & 2.07$\cdot 10^{9}$  & 5.16$\cdot 10^{-10}$ \\
\hline
Left common carotid artery (LCA) & 1.78$\cdot 10^{8}$ & 3.01$\cdot 10^{9}$ & 3.52$\cdot 10^{-10}$ \\
\hline
Left subclavian artery (LSA) & 7.09$\cdot 10^{7}$  & 1.19$\cdot 10^{9}$  & 9.35$\cdot 10^{-10}$ \\
\hline
Descending aorta (DA) & 7.8$\cdot 10^{6}$ & 1.31$\cdot 10^{8}$ & 7.72$\cdot 10^{-9}$ \\
\hline
\end{tabular}
\label{table_windkessel_aorta}
\end{table}

The computation of the Reynolds number (see equation \eqref{eq:Re}) is based on the diameter of the ascending aorta and on the velocity. 
Since the velocity is time-dependent (see  Figure~\ref{fig:BC_aorta}), $Re \in [0, 4200]$ over the cardiac cycle. 

We use an adaptive time step accordingly to the condition $CFL_{max} =  0.9$ where $CFL$ is the Courant-Friedrichs-Lewy number \cite{computing1969applied}. 
We run the simulations over $(t_0, T] = (0, 6] $ s and we observe that, after $8$ cardiac cycles (i.e., after 4.8 s), the system reaches a pseudo steady-state, indicating that the transient effects have been overcome. Consequentely, we consider the data associated to the final cardiac cycle, i.e. [5.4, 6] s. 
 At reduced order level, due to the better performance demonstrated for the idealized blood vessel (see Section~\ref{sec:cylinder}), the supremizer approach is adopted to stabilize the ROM.
Unlike what has been done for \emph{Case 1} here we compute the absolute reconstruction error for the velocity and the pressure 
\begin{equation}
    \bar{\varepsilon}_u = {\| \bm  u_{\text{FOM}}- \bm u_{\text{ROM}}\|_{L^2(\Omega)}}\qquad \text{and} \qquad \bar{\varepsilon}_p = {\| p_{\text{FOM}}-  p_{\text{ROM}}\|_{L^2(\Omega)}}.
    \label{eq:errors-absolute-aorta}
\end{equation}
\begin{figure}
	\centering
    \includegraphics[width=.4\textwidth]{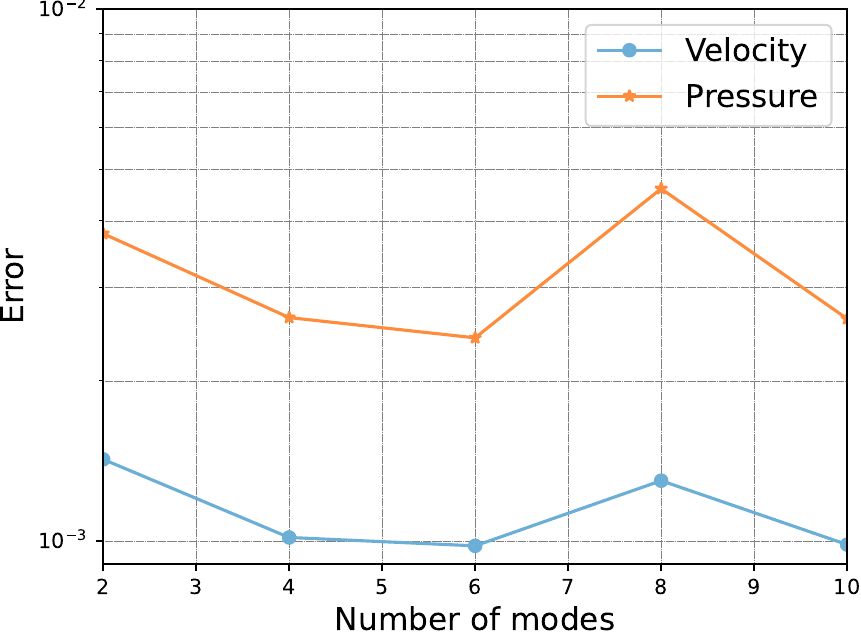}
	\caption{\textit{Case 2}: Time-averaged reconstruction error for velocity and pressure as the number of modes $N_{\Phi}$ increases.}
	\label{errup_N_aorta}
\end{figure}
This is because the relative error exhibits some peaks caused by small values of velocity and/or pressure at specific time instances. However, the outcomes presented demonstrate that the absolute error is a good metric for estabilishing the quality of the ROM reconstruction. 

We collect a total of 100 snapshots in time over [5.4, 6] s for training our ROM. For sake of clearness, we explicitly list the time steps collected in the training stage:
\begin{equation}
    \left\{5.4,5.41,5.42,\dots,5.98,5.99,6\right\} s.
    \label{time_set}
\end{equation}

\begin{figure}
	\centering
    \subfloat[][\label{fig:erru_aorta_proj}]{\includegraphics[width=.4\textwidth]{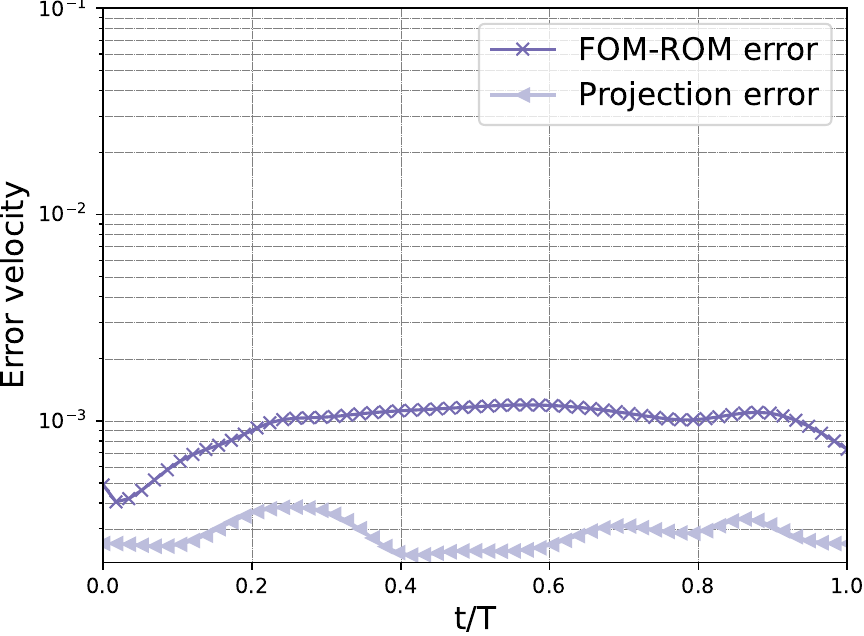}}\hspace{2ex}
 	\subfloat[][\label{fig:error_dtp_proj}]{\includegraphics[width=.4\textwidth]{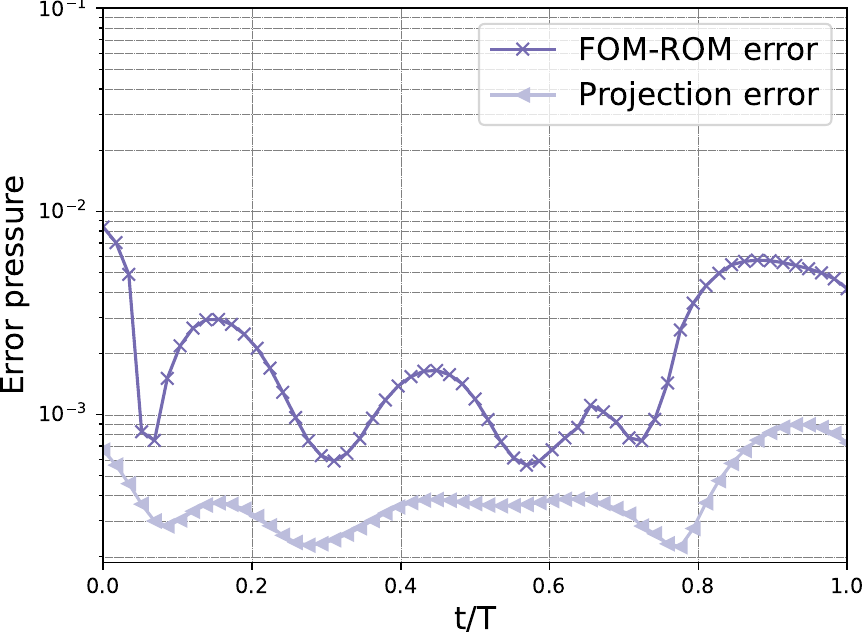}}\\
	\caption{\textit{Case 2}: Comparison between reconstruction error and projection error for velocity and pressure with $N_{\Phi} = 6$. The supremizer approach is adopted and $N_{\Phi}$ is the same for pressure, velocity and supremizers.}
	\label{errup_aorta_proj}
\end{figure}

As a first numerical experiment, we do not consider any new test point in the online phase, i.e. the validation set coincides with the training one \eqref{time_set}.  To reach $99.99\%$ of the cumulative energy of the eigenvalues, we need 12 modes for the velocity and only 1 mode for the pressure.
In Figure \ref{errup_N_aorta} we plot the time-averaged reconstruction error as $N_\Phi$ varies.  We set $N_\Phi = 6$ because the error reaches its minimum value.
The comparison between reconstruction error and projection error is reported in Figure~\ref{errup_aorta_proj}. 
For both velocity and pressure, the reconstruction error closely matches the projection error, demonstrating the capability of our ROM method in achieving an accurate approximation of the dynamics of the flow field. In particular, for the velocity, the reconstruction error is of the order $10^{-3}$ over the entire time window whilst the error for the pressure ranges between $10^{-3}$ and $10^{-2}$.

\begin{figure}
	\centering
    \subfloat[][\label{fig:erru_aorta_noliftP}]{\includegraphics[width=.4\textwidth]{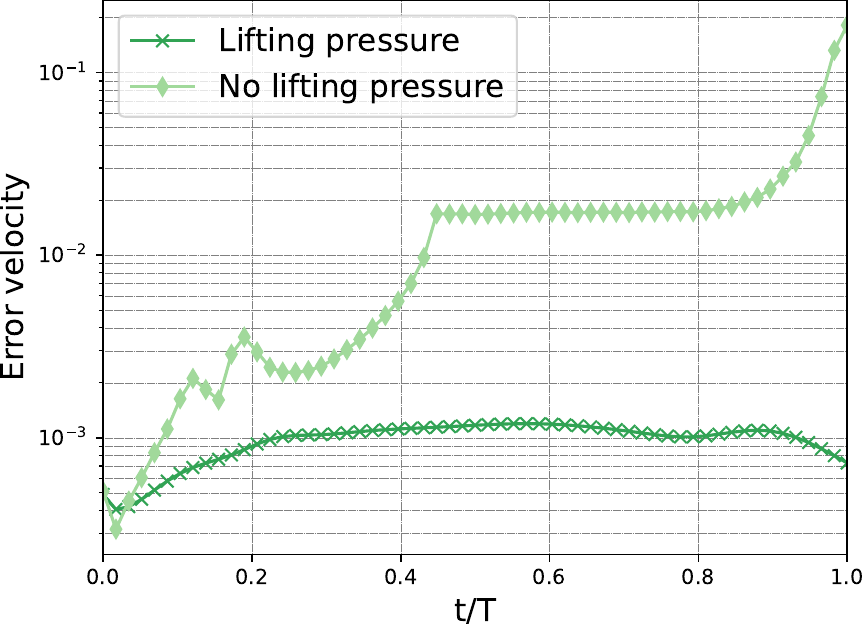}}
 	\hspace{2ex}\subfloat[][\label{fig:error_dtp_noliftP}]{\includegraphics[width=.4\textwidth]{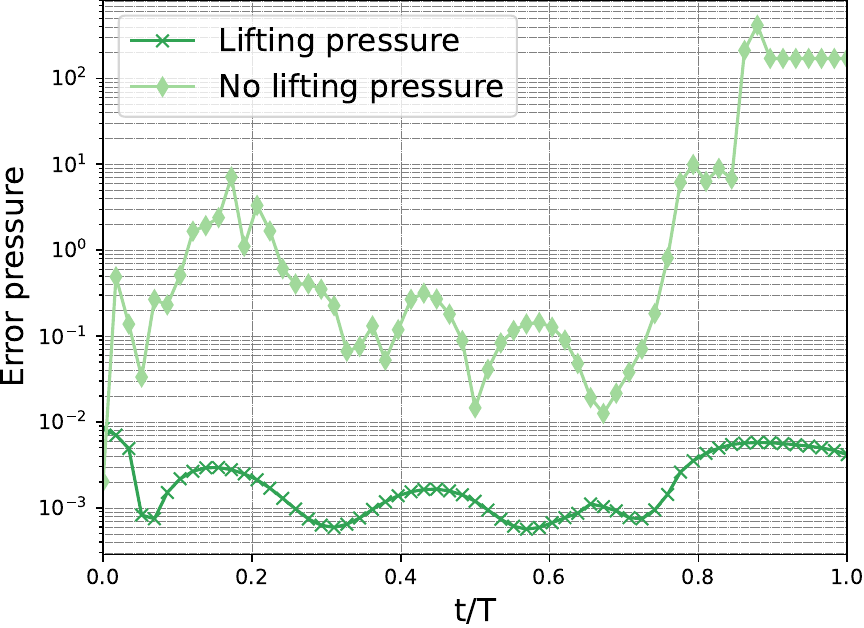}}\\
	\caption{\textit{Case 2}: Comparison of FOM-ROM errors for velocity and pressure with $N_{\Phi} = 6$, both with and without the use of the lifting function for the pressure.}
	\label{errup_aorta_noliftP}
\end{figure}

In Figure~\ref{errup_aorta_noliftP} we demonstrate the importance of incorporating the lifting function for the pressure in our ROM framework. In Figure~\ref{fig:error_dtp_noliftP} we note that the pressure error exceeds $10^2$ when the lifting function for the pressure is not employed. Due to the coupling between the velocity and pressure, the velocity error shown in Figure~\ref{fig:erru_aorta_noliftP}  also increases by two orders of magnitude, reaching $10^{-1}$.
It should be noted that we do not investigate the error computed without the lifting function for the velocity, as its crucial role has been thoroughly established in literature \cite{star2019novel}.

\begin{figure}[!htb]
    \centering
        \includegraphics[scale=0.28, trim = 0cm 0.85cm 0cm 0cm]{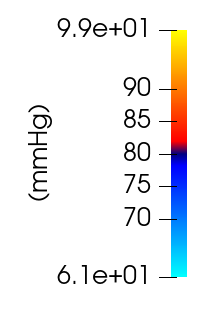}
        \begin{overpic}[width=0.2\textwidth]{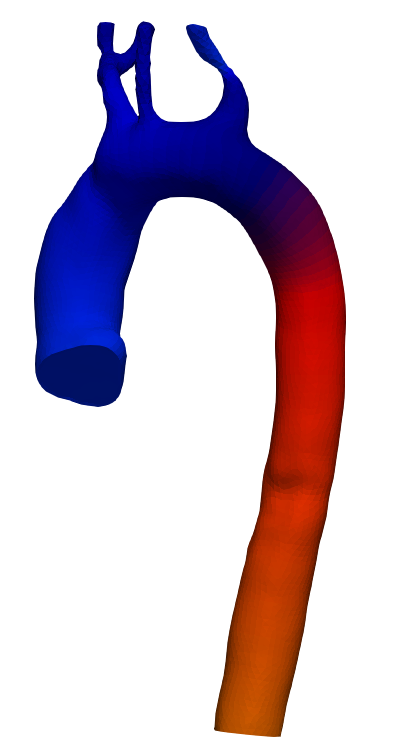} 
        \put(16,100){$p$ FOM}
        \end{overpic}
        \begin{overpic}[width=0.2\textwidth]{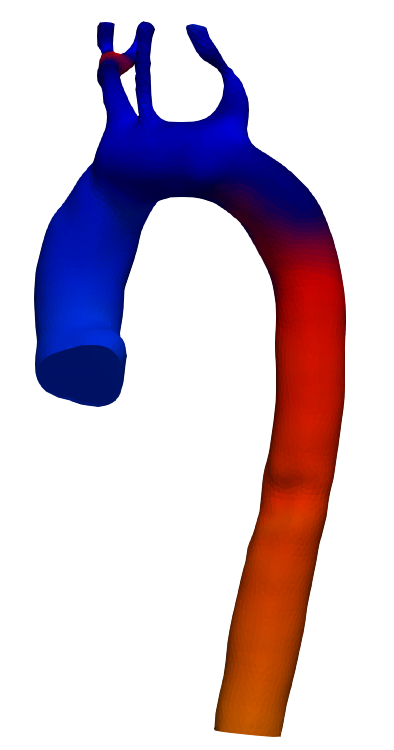} 
        \put(16,100){$p$ ROM}
        \put(15,109){$t=5.5$ s}
        \end{overpic}
        \vline
        \begin{overpic}[width=0.2\textwidth]{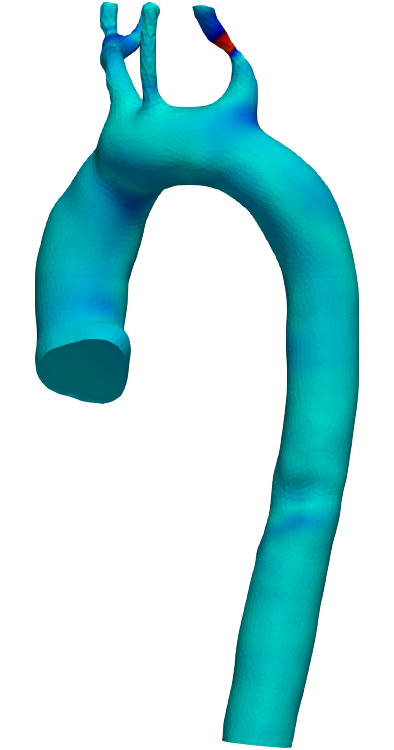} 
        \put(5,100){$\mid$ $p$ FOM - $p$ ROM $\mid$}
        \end{overpic}
        \hspace{2ex}
        \includegraphics[scale=0.28, trim = 0cm 0.85cm 0cm 0cm]{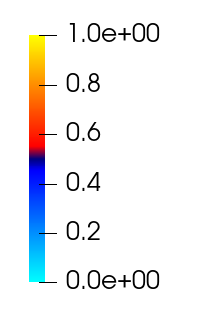}
        \\\vspace{8ex}
        \includegraphics[scale=0.28, trim = 0cm 0.85cm 0cm 0cm]{img/legend_p_mmHg.png}
        \begin{overpic}[width=0.2\textwidth]{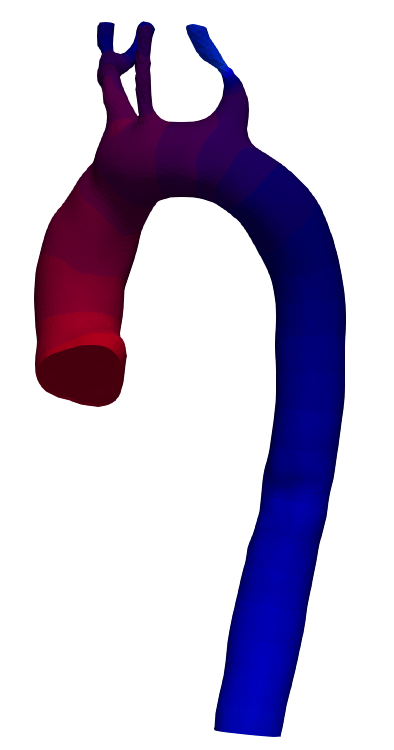} 
        \put(16,100){$p$ FOM}
        \end{overpic}
        \begin{overpic}[width=0.2\textwidth]{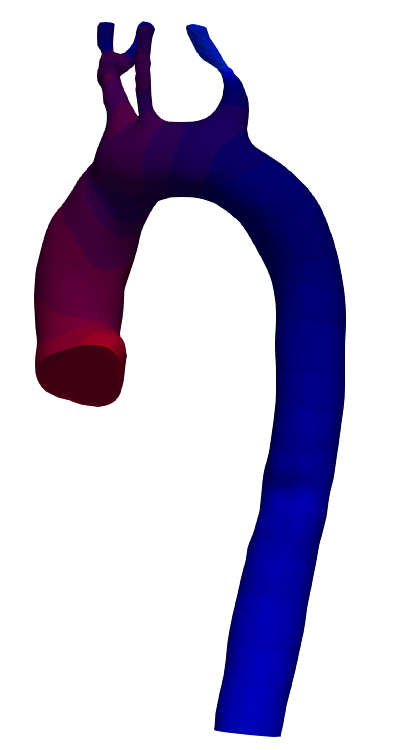} 
        \put(16,100){$p$ ROM}
        \put(15,109){$t=5.6$ s}
        \end{overpic}
        \vline
        \begin{overpic}[width=0.2\textwidth]{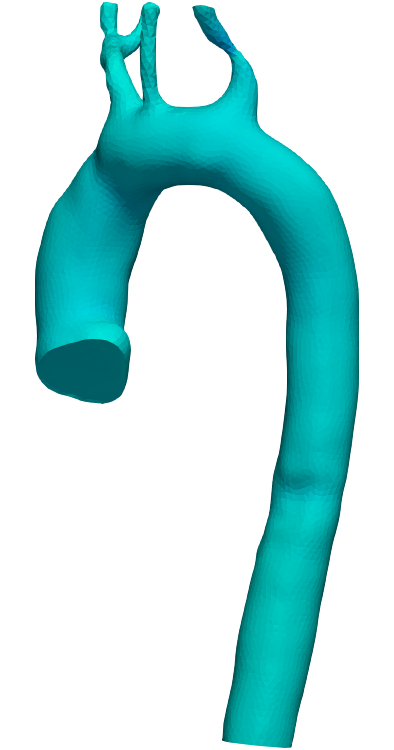} 
        \put(5,100){$\mid$ $p$ FOM - $p$ ROM $\mid$}
        \end{overpic}
        \hspace{2ex}
        \includegraphics[scale=0.28, trim = 0cm 0.85cm 0cm 0cm]{img/legend_err_p_aorta.png}
        \\\vspace{8ex}
        \includegraphics[scale=0.28, trim = 0cm 0.85cm 0cm 0cm]{img/legend_p_mmHg.png}
        \begin{overpic}[width=0.2\textwidth]{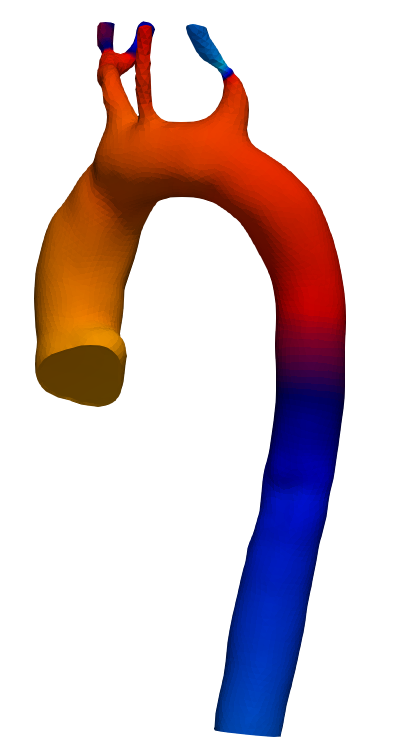} 
        \put(16,100){$p$ FOM}
        \end{overpic}
        \begin{overpic}[width=0.2\textwidth]{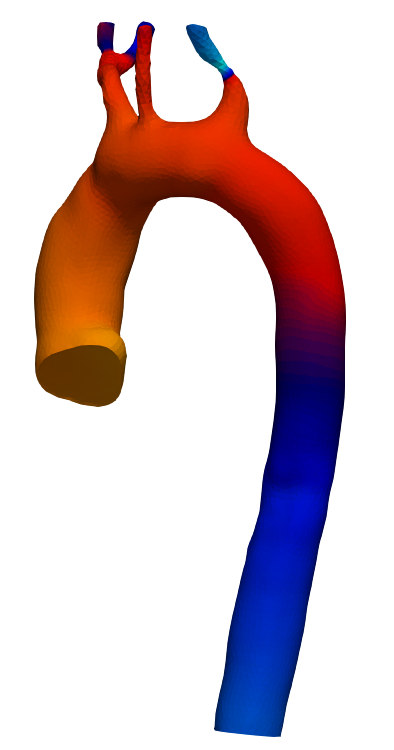} 
        \put(16,100){$p$ ROM}
        \put(15,109){$t=5.9$ s}
        \end{overpic}
        \vline
        \begin{overpic}[width=0.2\textwidth]{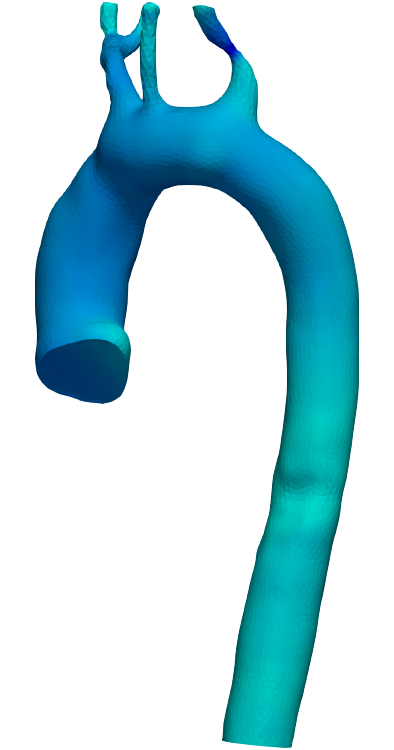} 
        \put(5,100){$\mid$ $p$ FOM - $p$ ROM $\mid$}
        \end{overpic}
        \hspace{2ex}
        \includegraphics[scale=0.28, trim = 0cm 0.85cm 0cm 0cm]{img/legend_err_p_aorta.png}
        
\caption{\textit{Case 2}: Qualitative comparison of FOM-ROM pressure at $t=5.5, 5.6, 5.9$ s.}
\label{fig:p-fom-rom-aorta}
\end{figure}

\begin{figure}[!htb]
    \centering
        \includegraphics[scale=0.28, trim = 0cm 0.85cm 0cm 0cm]{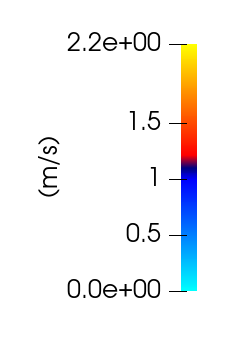}
        \begin{overpic}[width=0.2\textwidth]{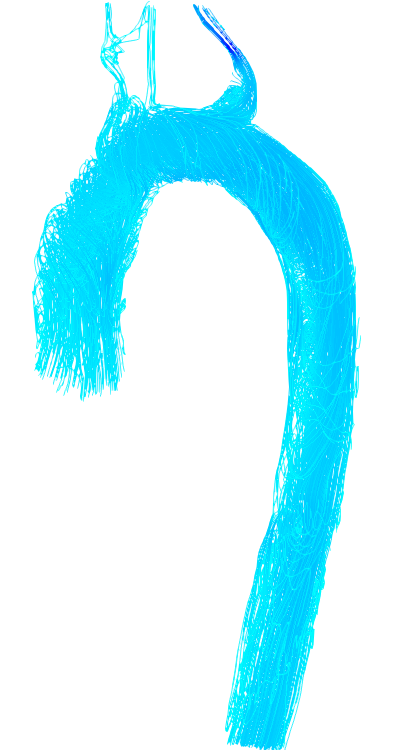} 
        \put(16,100){$\bm u$ FOM}
        \end{overpic}
        \begin{overpic}[width=0.2\textwidth]{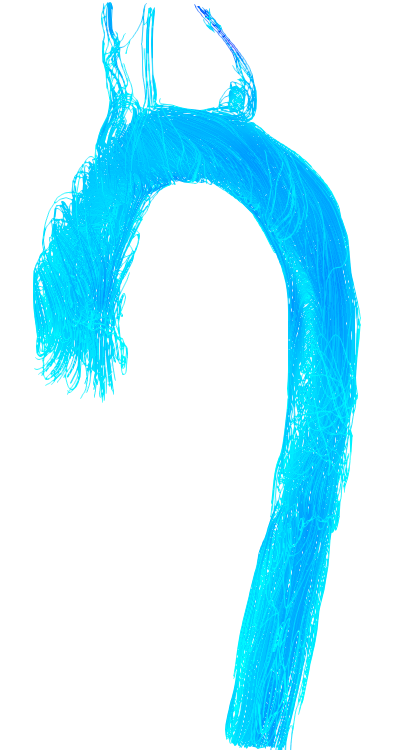} 
        \put(16,100){$\bm u$ ROM}
        \put(-10,109){$t=5.5$ s}
        \end{overpic}
        \\\vspace{8ex}
        \includegraphics[scale=0.28, trim = 0cm 0.85cm 0cm 0cm]{img/legend_u_reversed.png}
        \begin{overpic}[width=0.2\textwidth]{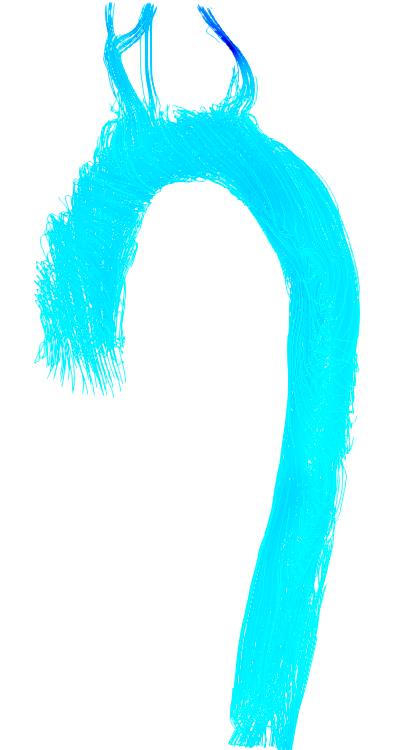} 
        \put(16,100){$\bm u$ FOM}
        \end{overpic}
        \begin{overpic}[width=0.2\textwidth]{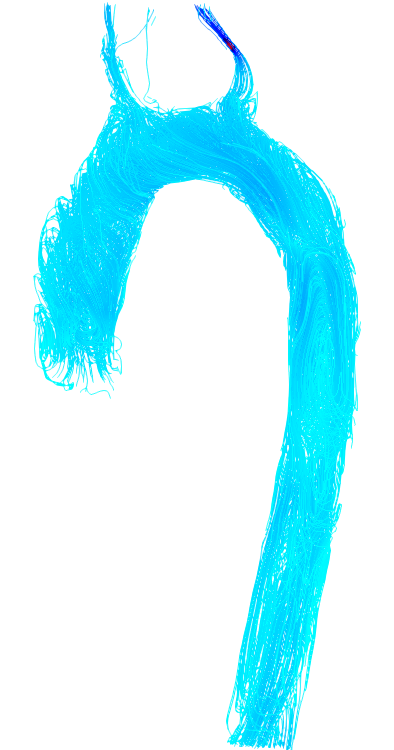} 
        \put(16,100){$\bm u$ ROM}
        \put(-10,109){$t=5.6$ s}
        \end{overpic}
        \\\vspace{8ex}
        \includegraphics[scale=0.28, trim = 0cm 0.85cm 0cm 0cm]{img/legend_u_reversed.png}
        \begin{overpic}[width=0.2\textwidth]{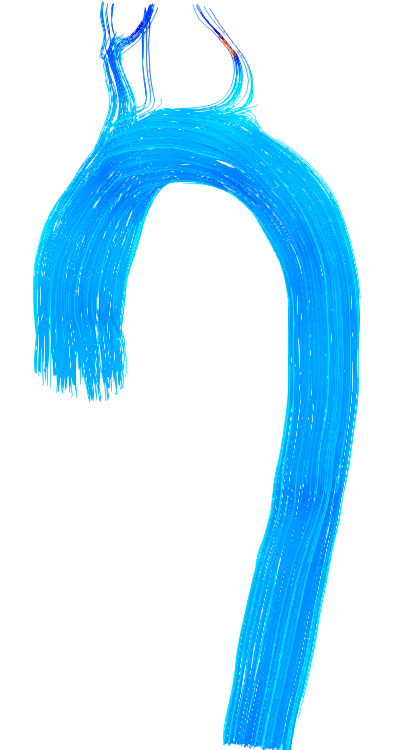} 
        \put(16,100){$\bm u$ FOM}
        \end{overpic}
        \begin{overpic}[width=0.2\textwidth]{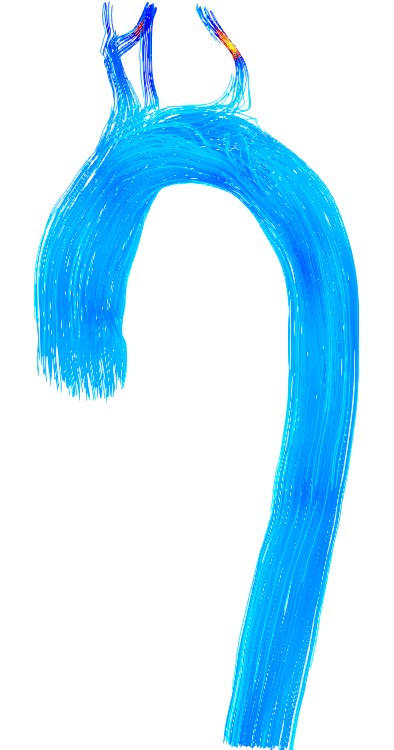} 
        \put(16,100){$\bm u$ ROM}
        \put(-10,109){$t=5.9$ s}
        \end{overpic}
\caption{\textit{Case 2}: Qualitative comparison of FOM-ROM velocity at $t=5.5, 5.6, 5.9$ s.}
\label{fig:u-fom-rom-aorta}
\end{figure}

\begin{figure}[!htb]
    \centering
        \includegraphics[width=.2\textwidth]{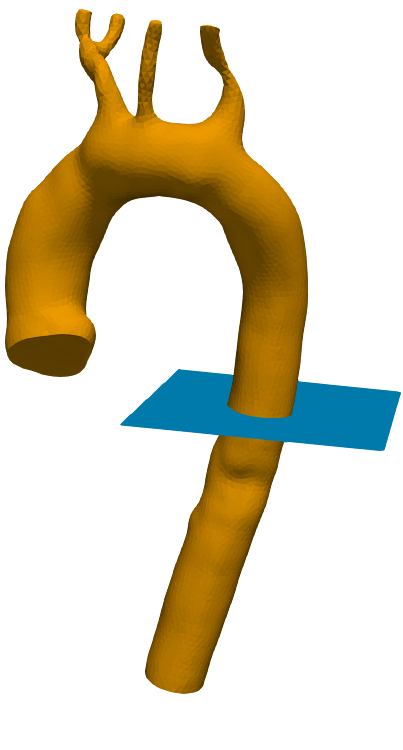} \\\vspace{1ex}
        \includegraphics[scale=0.26, trim = 0cm 0.85cm 0cm 0cm]{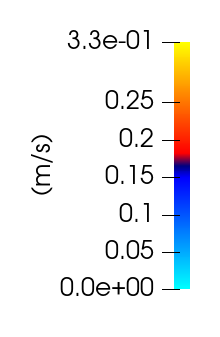}
        \begin{overpic}[width=0.18\textwidth]{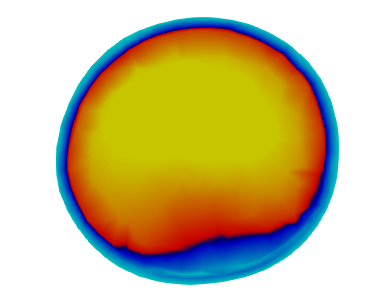} 
        \put(29,85){$\bm u$ FOM}
        \end{overpic}
        \begin{overpic}[width=0.18\textwidth]{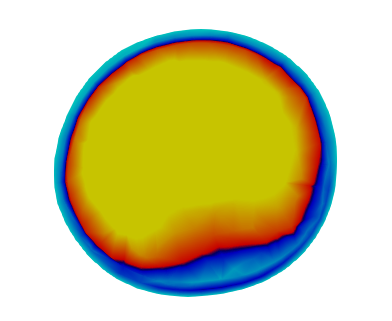} 
        \put(29,86){$\bm u$ ROM}
        \put(28,109){$t=5.5$ s}
        \end{overpic}
        \vline
        \hspace{4ex}
        \begin{overpic}[width=0.155\textwidth]{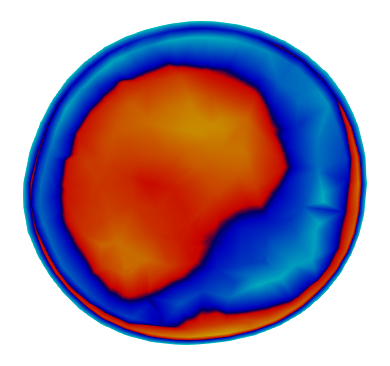} 
        \put(-10,100){$\mid$ $\bm u$ FOM - $\bm u$ ROM $\mid$}
        \end{overpic}
        \hspace{4ex}
        \includegraphics[scale=0.28, trim = 0cm 0.85cm 0cm 0cm]{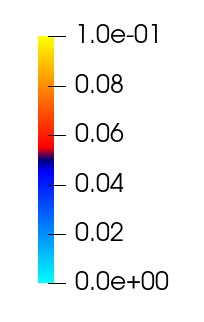}
        \\\vspace{8ex}
        \includegraphics[scale=0.26, trim = 0cm 0.85cm 0cm 0cm]{img/legend_slice_u.png}
        \begin{overpic}[width=0.18\textwidth]{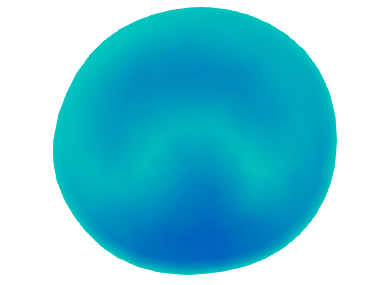} 
        \put(29,85){$\bm u$ FOM}
        \end{overpic}
        \begin{overpic}[width=0.18\textwidth]{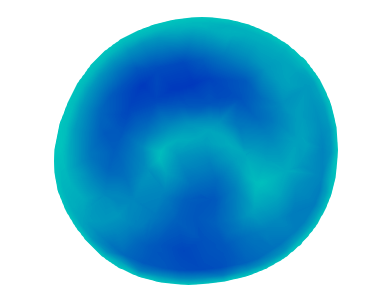} 
        \put(29,85){$\bm u$ ROM}
        \put(28,109){$t=5.6$ s}
        \end{overpic}
        \vline
        \hspace{4ex}
        \begin{overpic}[width=0.153\textwidth]{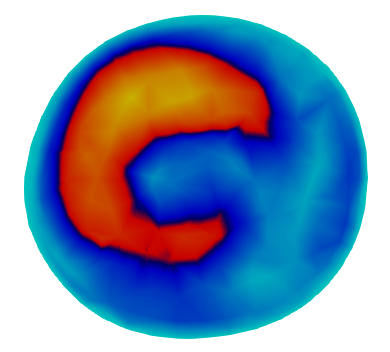} 
        \put(-10,100){$\mid$ $\bm u$ FOM - $\bm u$ ROM $\mid$}
        \end{overpic}
        \hspace{4ex}
        \includegraphics[scale=0.28, trim = 0cm 0.85cm 0cm 0cm]{img/legend_error_slice_u.png}
        \\\vspace{8ex}
        \includegraphics[scale=0.26, trim = 0cm 0.85cm 0cm 0cm]{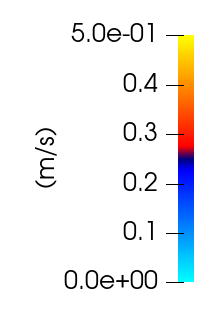}
        \begin{overpic}[width=0.18\textwidth]{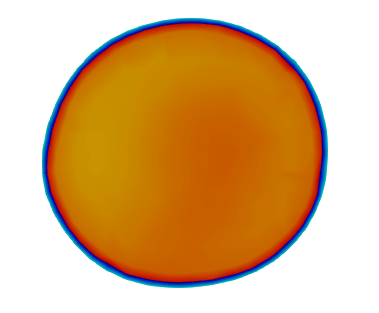} 
        \put(29,87){$\bm u$ FOM}
        \end{overpic}
        \begin{overpic}[width=0.18\textwidth]{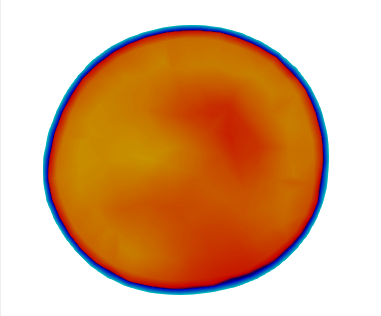} 
        \put(29,87){$\bm u$ ROM}
        \put(28,109){$t=5.9$ s}
        \end{overpic}
        \vline
        \hspace{4ex}
        \begin{overpic}[width=0.16\textwidth]{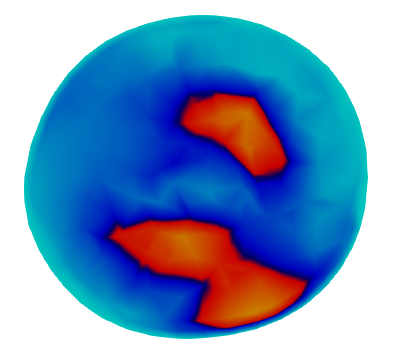} 
        \put(-10,100){$\mid$ $\bm u$ FOM - $\bm u$ ROM $\mid$}
        \end{overpic}
        \hspace{4ex}
        \includegraphics[scale=0.28, trim = 0cm 0.85cm 0cm 0cm]{img/legend_error_slice_u.png}
        
\caption{\textit{Case 2}: Qualitative comparison of FOM-ROM velocity at $t=5.5, 5.6, 5.9$ s  on the slice (blue plane) shown in the descending aorta.}
\label{fig:u-fom-rom-aorta_slice}
\end{figure}

Qualitative comparisons between FOM and ROM are displayed in Figure~\ref{fig:p-fom-rom-aorta}, \ref{fig:u-fom-rom-aorta} and \ref{fig:u-fom-rom-aorta_slice} at $t=5.5, 5.6, 5.9$ s (i.e. for training data). Moreover, we also show the spatial distribution of the error. In Figure~\ref{fig:p-fom-rom-aorta} we appreciate a very good performance of our ROM in the reconstruction of the pressure. 
In Figure~\ref{fig:u-fom-rom-aorta} the FOM and ROM streamlines for the velocity field are depicted.
We appreciate a good matching between the two solutions.
For further comparison in terms of velocity,  Figure~\ref{fig:u-fom-rom-aorta_slice} shows the FOM and ROM velocity magnitude on a slice of the descending aorta. 
We see that the ROM effectively captures the flow structures with small discrepancies compared to the FOM.


So far we have considered only the training data for the validation of our ROM, i.e. without adding any new time instance in the online phase. Now we are going to perform this kind of analysis. The set defined in \eqref{time_set} continues to be employed during the offline phase, whereas the following set is used for the evaluation stage:
\begin{equation}
    \{5.4, 5.405, 5.41, 5.415,\dots, 5.985, 5.9, 5.995, 6  \}.
    \label{new_time}
\end{equation}

\begin{figure}
	\centering
    \includegraphics[width=.4\textwidth]{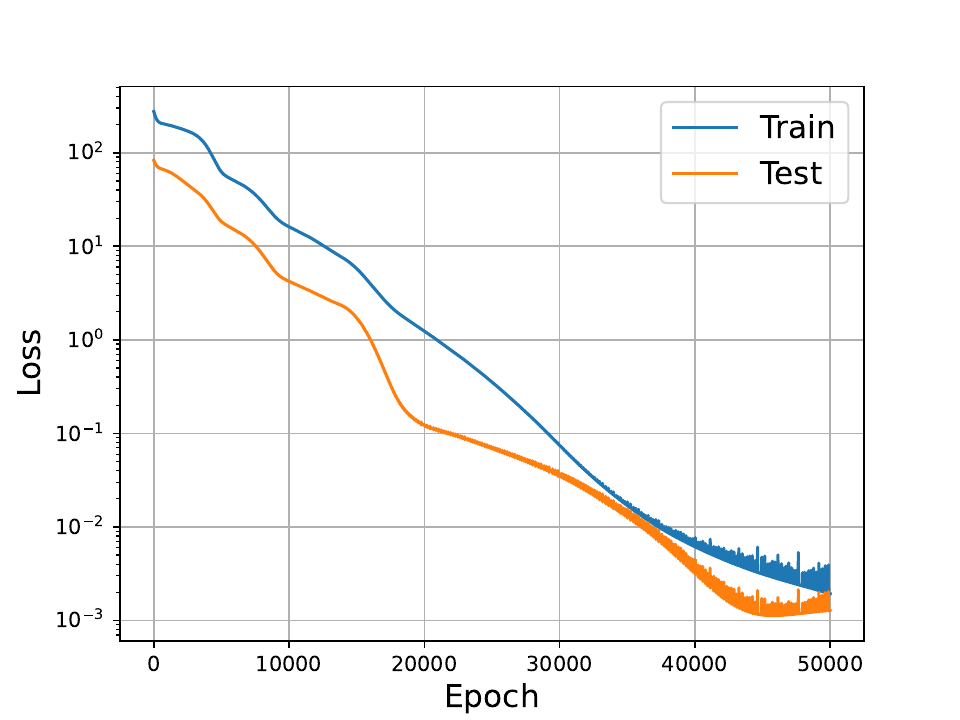}
	\caption{\textit{Case 2}: Loss function of the neural network computed for the train and for the test set.}
	\label{loss}
\end{figure}

Table \ref{tab:nn} reports the hyperparameters of the neural network used to interpolate the outflow pressure in the new time instances. The number of hidden layers is set equal to 2, as it provides a good balance between accuracy and efficiency, allowing the network to capture intricate patterns in the data without becoming overly complex \cite{hornik1989multilayer}. The tuning of the other hyperparameters is performed as follows: given a specific activation function and number of hidden layers, we increase the values of the learning rate and the number of hidden neurons until a satisfactory decay of the loss function is achieved and overfitting is avoided. The dataset \eqref{time_set} is normalized and split into train ($80\%$ of the dataset) and test ($20\%$ of the dataset) sets.
Figure~\ref{loss} illustrates the decrease of the loss function during the training process for both datasets. We observe that the minimization of the loss function is similar for both the training and test set, indicating that the model is generalizing well. This suggests that it is not overfitting to the training data and is capable of making accurate predictions on new time steps. 
\begin{table}
\centering
\renewcommand{\arraystretch}{1.5}
\caption{\textit{Case 2}: Hyperparameters used for the feedforward neural network.}
\begin{tabular}{ccccc}
\hline
\rowcolor{gray!20} Neurons per layer & Activation funtion  & Number of epochs & Learning rate & Hidden layers  \\
\hline
150  &   Softplus   & 50000 & 5$\cdot 10^{-6}$ & 2
\\
\hline
\end{tabular}
\label{tab:nn}
\end{table}

\begin{figure}
	\centering
    \subfloat[][\label{fig:erru_aorta_newt}]{\includegraphics[width=.4\textwidth]{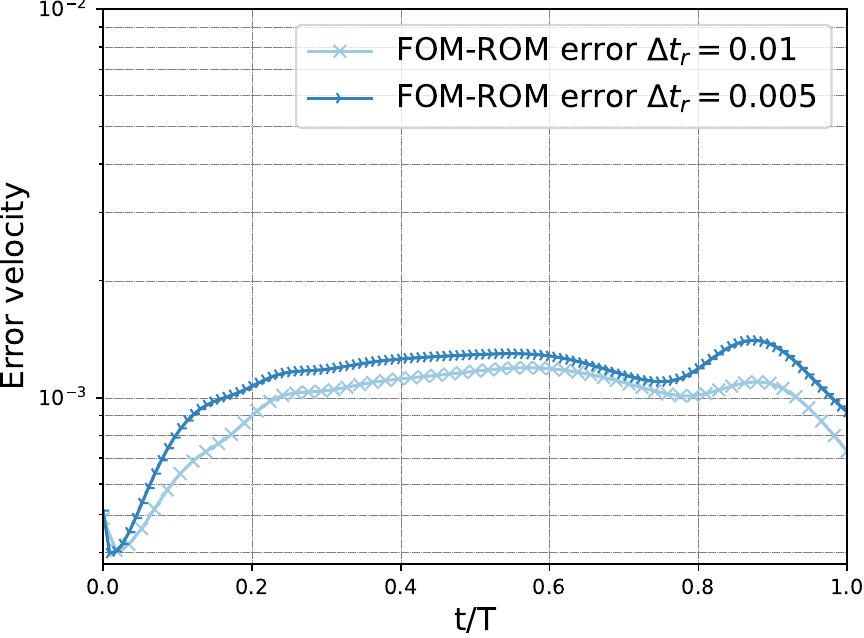}}\hspace{2ex}
 	\subfloat[][\label{fig:error_dtp_newt}]{\includegraphics[width=.4\textwidth]{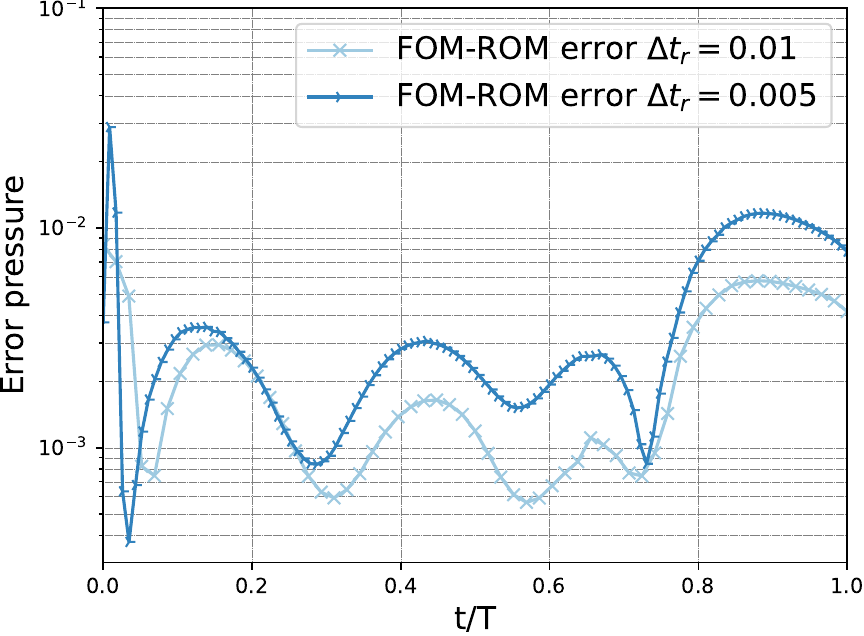}}\\
	\caption{\textit{Case 2}: Comparison of FOM-ROM errors for velocity and pressure with $N_{\Phi} = 6$ for the two numerical experiments performed}
    \label{errup_aorta_newt}
\end{figure}

In Figure~\ref{errup_aorta_newt}, we compare the  FOM-ROM error 
for velocity and pressure evaluating the  reduced solution for $\Delta t_r = 0.01$ (i.e. the former numerical experiment described in this section)  and $ \Delta t_r = 0.005$ (i.e. the latter numerical experiment described in this section). We note a slight increase in the velocity and pressure error for the test data (Figure~\ref{fig:erru_aorta_newt}). Nevertheless, the error trends exhibit strong similarity, affirming the robustness of our ROM framework, which remains applicable for predicting time instances where the outflow pressure at full order level has not been stored. 

\begin{figure}[!htb]
    \centering
    \vspace{3ex}
        \includegraphics[scale=0.28, trim = 0cm 0.85cm 0cm 0cm]{img/legend_p_mmHg.png}
        \begin{overpic}[width=0.2\textwidth]{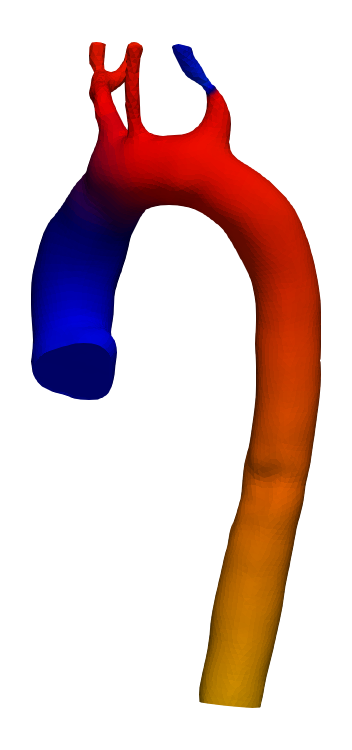} 
        \put(16,100){$p$ FOM}
        \end{overpic}
        \begin{overpic}[width=0.2\textwidth]{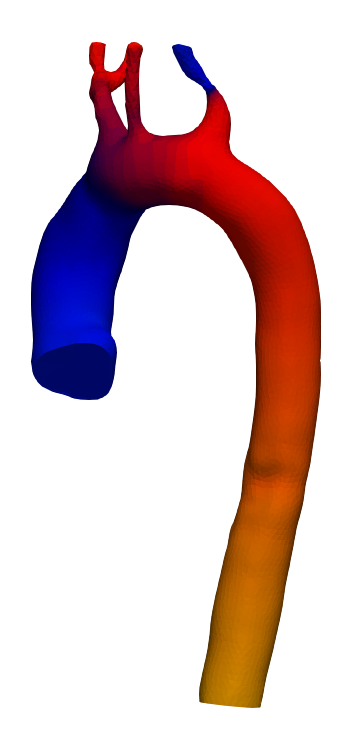} 
        \put(16,100){$p$ ROM}
        \end{overpic}
        \vline
        \begin{overpic}[width=0.2\textwidth]{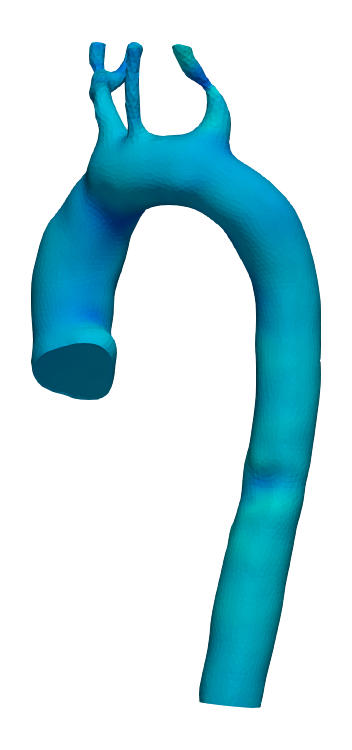} 
        \put(5,100){$\mid$ $p$ FOM - $p$ ROM $\mid$}
        \end{overpic}
        \hspace{2ex}
        \includegraphics[scale=0.28, trim = 0cm 0.85cm 0cm 0cm]{img/legend_err_p_aorta.png}
            
\caption{\textit{Case 2}: Qualitative comparison of FOM-ROM pressure for the new time $t=5.475$ s.}
\label{fig:p-fom-rom-aorta-newt}
\end{figure}

\begin{figure}[!htb]
    \centering
        \includegraphics[scale=0.28, trim = 0cm 0.85cm 0cm 0cm]{img/legend_u_reversed.png}
        \begin{overpic}[width=0.2\textwidth]{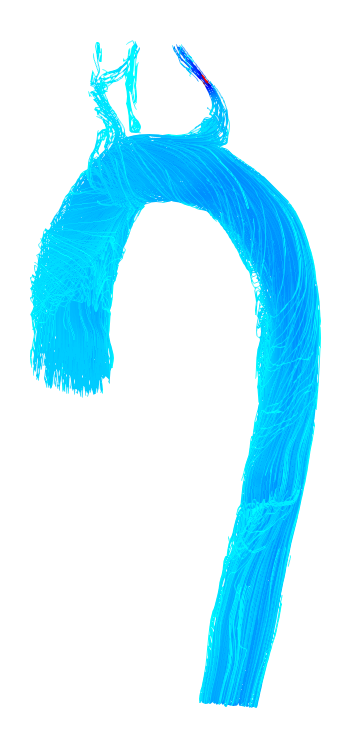} 
        \put(16,100){$\bm u$ FOM}
        \end{overpic}
        \begin{overpic}[width=0.2\textwidth]{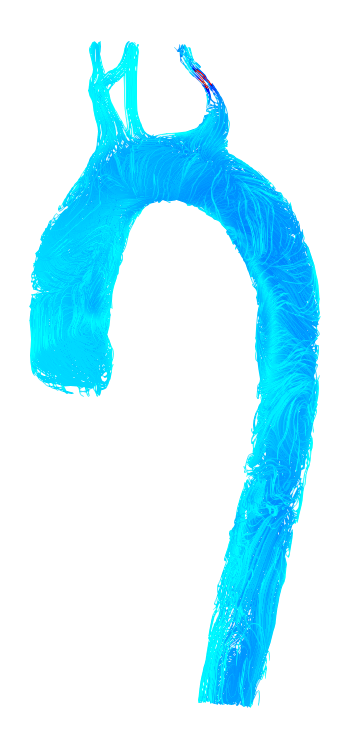} 
        \put(16,100){$\bm u$ ROM}
        \end{overpic}
\caption{\textit{Case 2}: Qualitative comparison of FOM-ROM velocity for the new time $t=5.475$ s.}
\label{fig:u-fom-rom-aorta-newt}
\end{figure}

\begin{figure}[!htb]
    \centering
        \includegraphics[scale=0.26, trim = 0cm 0.85cm 0cm 0cm]{img/legend_slice_u.png}
        \begin{overpic}[width=0.18\textwidth]{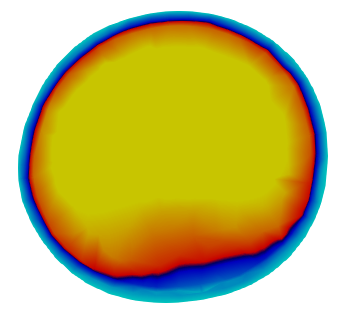} 
        \put(29,95){$\bm u$ FOM}
        \end{overpic}
        \begin{overpic}[width=0.18\textwidth]{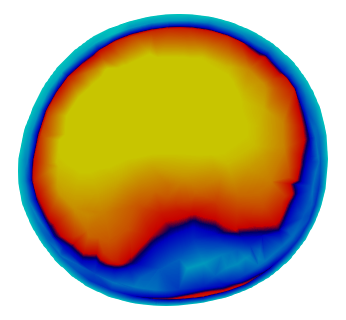} 
        \put(29,95){$\bm u$ ROM}
        \end{overpic}
        \vline
        \hspace{4ex}
        \begin{overpic}[width=0.18\textwidth]{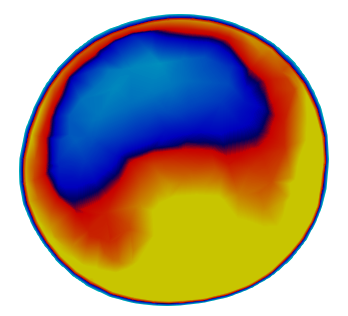} 
        \put(-10,95){$\mid$ $\bm u$ FOM - $\bm u$ ROM $\mid$}
        \end{overpic}
        \hspace{4ex}
        \includegraphics[scale=0.28, trim = 0cm 0.85cm 0cm 0cm]{img/legend_error_slice_u.png}
\caption{\textit{Case 2}: Qualitative comparison of FOM-ROM velocity for the new $t=5.475$ s.}
\label{fig:u-fom-rom-aorta_slice-newt}
\end{figure}

For the sake of completeness we also provide qualitative comparisons at a new time step $t = 5.475$ s. In Figure~\ref{fig:p-fom-rom-aorta-newt} we notice a good agreement in the pressure distribution between the FOM and the ROM. 
Figure~\ref{fig:u-fom-rom-aorta-newt} shows the FOM and ROM streamlines of velocity and Figure~\ref{fig:u-fom-rom-aorta_slice-newt} the internal velocity, on the same slice introduced in Figure~\ref{fig:u-fom-rom-aorta_slice}. We observe that the FOM-ROM velocity distributions are quite similar; however, the reconstruction is slightly worse for this new time step, which is reflected in the higher error obtained for new data in Figure~\ref{fig:p-fom-rom-aorta-newt} for $0<t/T<0.2$.

Finally, we comment on the computational cost using an Intel(R) Core(TM) i7-7700 CPU @ 3.60GHz 16GB RAM. Each FOM simulation requires about
1 day, $27$ seconds for the POD algorithm to complete and $87$ seconds for the training of the neural network, while less than 1 second is needed for the evaluation phase (285 milliseconds). Therefore, the speedup, which is defined as the ratio between the CPU time taken by the FOM simulation and the CPU time taken by the solution of the dynamic system for the reduced coefficients, is significant (in the order of $\mathcal{O}(10^{5})$).

\section{Conclusions and perspectives}\label{sec:conclusion}
This paper treats an hybrid model order reduction technique for cardiovascular application. The classic POD-Galerkin procedure is used in combination with data-driven techniques to generalize our results. Nonhomogeneous boundary conditions, which are not automatically preserved at ROM level, are treated with the lifting function method \cite{saddam2017pod,graham1999optimal,gunzburger2007reduced}. This methodology shift the nonhomogeneous high fidelity solutions to obtain an homogeneous reduced space. The boundary conditions are added during the reconstruction phase, so that every boundary value can be integrated in the ROM framework. In the context of CFD and Navier-Stokes equation for real-world problems, the lifting function method is well-known for the velocity field. 
Velocity profiles are typically imposed with an available functional form. However, this method remains largely unexplored for pressure, which often assumes homogeneous boundary conditions. For realistic simulations, many cardiovascular applications require non-homogeneous pressure conditions. At the FOM level, the 3D Windkessel model is frequently used to impose outflow pressure in vessels. This work introduce the outflow pressure obtained from Windkessel in the POD-Galerkin framework, with the lifting function approach. 
Since we do not have a function for the pressure, neural networks are employed to interpolate it at the outlet once the FOM is trained. This approach allows the ROM to provide solutions for new parameter values, beyond just the data available from the FOM.

The framework resumed here is tested on an idealized benchmark (a cylinder) and on a 3D patient-specific aortic arch. The accuracy of our ROM is tested in terms of errors and qualitative comparisons of pressure and velocity. The results are promising and the qualitative comparison are very similar. Also the reduction in term of time is good, because we need less than one second to obtain the ROM solution, once the FOM is trained.

However, significant effort is still required to advance ROMs for cardiovascular applications. Specifically, incorporating machine-learning techniques can be highly beneficial for real patient-specific studies. 
Future research could explore the introduction of autoencoders \cite{phillips2021autoencoder}, which provide a nonlinear alternative to POD and have the potential to more efficiently capture features or patterns in high-fidelity model results. 
Another area of future research could involve studying blood flow using deep neural networks, particularly physics-informed neural networks, incorporating both physical and geometric parameterizations. While the results presented in this study provide valuable insights, it is important to note that using a single patient as a test case represents a limitation of our analysis, highlighting the need of further investigations involving a larger set to enhance the robustness and applicability of our results.

\section*{Acknowledgments}
We acknowledge the support provided by PRIN “FaReX - Full and Reduced order modeling of
coupled systems: focus on non-matching methods and automatic learning” project, PNRR NGE iNEST “Interconnected Nord-Est Innovation Ecosystem” project, INdAM-GNCS 2019–2020 projects
and PON “Research and Innovation on Green related issues” FSE REACT-EU 2021 project.





\printbibliography
 
\end{document}